\documentclass{article}
\usepackage{amsmath}
\usepackage{t1enc}
\usepackage[latin1]{inputenc}
\usepackage[english]{babel}
\usepackage{amsmath,amsthm}
\usepackage{amsfonts}
\usepackage{latexsym}
\usepackage[dvips]{graphicx}
\usepackage{graphicx}
\usepackage[natural]{xcolor}

\newtheorem{theorem}{Theorem}
\newtheorem{lemma}[theorem]{Lemma}
\newtheorem{proposition}[theorem]{Proposition}

\newtheorem{definition}[theorem]{Definition}

\newtheorem{remark}[theorem]{Remark}

\input{tcilatex}
\begin{document}

\date{30/10/ 2012}
\author{Liangquan Zhang$^{1,2}$ \thanks{
This work was supported by Marie Curie Initial Training Network (ITN)
project: \textquotedblright Deterministic and Stochastic Controlled System
and Application\textquotedblright , FP7-PEOPLE-2007-1-1-ITN, No. 213841-2
and National Natural Science Foundation of China Grant 10771122, Natural
Science Foundation of Shandong Province of China Grant Y2006A08 and National
Basic Research Program of China (973 Program, No. 2007CB814900) and .
Corresponding author, E-mail: xiaoquan51011@163.com. } Yufeng Shi$^{1}$ \\
1. School of Mathematics, Shandong University\\
Jinan 250100, People's Republic of China.\\
{2. }Laboratoire de Math\'ematiques, \\
Universit\'e de Bretagne Occidentale, 29285 Brest C\'edex, France.}
\title{General Doubly Stochastic Maximum Principle and Its Applications to
Optimal Control of SPDEs}
\maketitle

\begin{abstract}
In this paper, we prove the necessary and sufficient maximum principles
(NSMPs in short) for the optimal control of systems described by a
quasilinear stochastic heat equation within convex control domains, which
all the coefficients contain control variables. For that, the optimal
control problem of fully coupled forward-backward doubly stochastic system
is studied. We apply our NSMPs to treat a kind of forward-backward doubly
stochastic linear quadratic optimal control problems and an example of
optimal control of stochastic partial differential equations (SPDEs in
short) as well.
\end{abstract}

\section{Introduction}

In order to provide a probabilistic interpretation for the solutions of a
class of quasilinear stochastic partial differential equations (SPDEs in
short), Pardoux and Peng [15] introduced the following backward doubly
stochastic differential equation (BDSDE in short): 
\begin{equation*}
Y_{t}=\xi +\int_{t}^{T}f(s,Y_{s},Z_{s})\text{d}s+\int_{t}^{T}g(s,Y_{s},Z_{s})%
\text{d}\overleftarrow{B_{s}}-\int_{t}^{T}Z_{s}\text{d}\overrightarrow{W_{s}}%
,\quad 0\leq t\leq T.
\end{equation*}%
Note that the integral with respect to $\{B_{t}\}$ is a \textquotedblleft
backward It\^{o} integral\textquotedblright\ and the integral with respect
to $\{W_{t}\}$ is a standard forward It\^{o} integral. These two types of
integrals are particular cases of the It\^{o}-Skorohod integral (for more
details see [11] and [15]). Pardoux and Peng [15] have obtained the
relationship between BDSDEs and a certain quasilinear stochastic partial
differential equations (SPDEs in short). More precisely 
\begin{equation*}
\left\{ 
\begin{array}{c}
u\left( t,x\right) =\varphi \left( x\right) +\int_{t}^{T}\left[ \mathcal{L}%
u\left( s,x\right) +f\left( s,x,u\left( s,x\right) ,\left( \nabla u\sigma
\right) \left( s,x\right) \right) \right] \text{d}s \\ 
+\int_{t}^{T}g\left( s,x,u\left( s,x\right) ,\left( \nabla u\sigma \right)
\left( s,x\right) \right) \text{d}\overleftarrow{B_{s}},\quad 0\leq t\leq T,%
\end{array}%
\right. 
\end{equation*}%
where $u:\left[ 0,T\right] \times \mathbb{R}^{d}\mathbb{\rightarrow R}^{k}$
where $d,$ $k\in N$, and $\nabla u\left( s,x\right) $ denotes the first
order derivative of $u\left( s,x\right) $ with respect to $x$, and 
\begin{equation*}
\mathcal{L}u=\left( 
\begin{array}{c}
Lu_{1} \\ 
\vdots  \\ 
Lu_{k}%
\end{array}%
\right) ,
\end{equation*}%
with 
\begin{equation*}
L\phi \left( x\right) =\frac{1}{2}\sum_{i,j=1}^{d}\left( \sigma \sigma
^{\ast }\right) _{ij}\left( x\right) \frac{\partial ^{2}\phi \left( x\right) 
}{\partial x_{i}\partial x_{j}}+\sum_{i=1}^{d}b_{i}\left( x\right) \frac{%
\partial \phi \left( x\right) }{\partial x_{i}}
\end{equation*}%
(for more details see [15]).

In 2003, Peng and Shi [18] introduced a type of time-symmetric
forward-backward stochastic differential equations, i.e., so-called fully
coupled forward-backward doubly stochastic differential equations (FBDSDE in
short): 
\begin{equation}
\left\{ 
\begin{array}{lll}
y_{t} & = & x+\int_{0}^{t}f\left( s,y_{s},Y_{s},z_{s},Z_{s}\right) \text{d}%
s+\int_{0}^{t}g\left( s,y_{s},Y_{s},z_{s},Z_{s}\right) \text{d}%
\overrightarrow{W_{s}}-\int_{0}^{t}z_{s}\text{d}\overleftarrow{B_{s}}, \\ 
Y_{t} & = & \varphi \left( y_{T}\right) +\int_{t}^{T}F\left(
s,y_{s},Y_{s},z_{s},Z_{s}\right) \text{d}s+\int_{t}^{T}G\left(
s,y_{s},Y_{s},z_{s},Z_{s}\right) \text{d}\overleftarrow{B_{s}}%
+\int_{t}^{T}Z_{s}\text{d}\overrightarrow{W_{s}}.%
\end{array}%
\right.   \tag{1.2}
\end{equation}

In FBDSDEs (1.2), the forward equation is \textquotedblleft
forward\textquotedblright\ with respect to a standard stochastic integral d$%
\overrightarrow{W_{t}}$, as well as \textquotedblleft
backward\textquotedblright\ with respect to a backward stochastic integral d$%
\overleftarrow{B_{t}}$; the coupled \textquotedblleft backward
equation\textquotedblright\ is \textquotedblleft forward\textquotedblright\
under the backward stochastic integral d$\overleftarrow{B_{t}}$ and
\textquotedblleft backward\textquotedblright\ under the forward one. In
other words, both the forward equation and the backward one are types of
BDSDE (1.1) with different directions of stochastic integrals. So (1.2)
provides a very general framework of fully coupled forward-backward
stochastic systems. Peng and Shi [18] proved the existence and uniqueness of
solutions to FBDSDE (1.2) with arbitrarily fixed time duration under some
monotone assumptions. FBDSDE (1.2) can provide a probabilistic
interpretation for the solutions of a general class of quasilinear SPDEs.

In this paper, we consider the following quasilinear SPDEs with control
variable: 
\begin{equation}
\left\{ 
\begin{array}{c}
u\left( t,x\right) =\varphi \left( x\right) +\int_{t}^{T}\left[ \mathcal{L}%
^{v}u\left( s,x\right) +F\left( s,x,u\left( s,x\right) ,\left( \nabla
u\sigma \right) \left( s,x,u\right) ,v\left( s\right) \right) \right] \text{d%
}s \\ 
+\int_{t}^{T}G\left( s,x,u\left( s,x\right) ,\left( \nabla u\sigma \right)
\left( s,x,u\right) ,v\left( s\right) \right) \text{d}\overleftarrow{B_{s}}%
,\quad 0\leq t\leq T,%
\end{array}%
\right.  \tag{1.3}
\end{equation}%
where $u:\left[ 0,T\right] \times \mathbb{R}^{d}\mathbb{\rightarrow R}^{k}$
and $\nabla u\left( s,x\right) $ denotes the first order derivative of $%
u\left( s,x\right) $ with respect to $x$, and 
\begin{equation*}
\mathcal{L}^{v}u=\left( 
\begin{array}{c}
L^{v}u_{1} \\ 
\vdots \\ 
L^{v}u_{k}%
\end{array}%
\right) ,
\end{equation*}%
with 
\begin{equation*}
L^{v}\phi \left( x\right) =\frac{1}{2}\sum_{i,j=1}^{d}\left( gg^{\ast
}\right) _{ij}\left( x,v\right) \frac{\partial ^{2}\phi \left( x\right) }{%
\partial x_{i}\partial x_{j}}+\sum_{i=1}^{d}f_{i}\left( x,v\right) \frac{%
\partial \phi \left( x\right) }{\partial x_{i}}.
\end{equation*}%
It is worth to pointing out that all the coefficients contain the control
variable. (For more details see Section 5).

Let us describe the problem solved in this paper. Set $\mathcal{U}_{ad}$ be
an admissible control set. The definitions of notations used here can be
found in Section 2. The optimal control problem of SPDEs (1.3) is to find an
optimal control $v^{\ast }\left( \cdot \right) \in \mathcal{U}_{ad}$, such
that 
\begin{equation*}
J\left( v^{\ast }\left( \cdot \right) \right) \doteq \underset{v\left( \cdot
\right) \in \mathcal{U}_{ad}}{\inf }J\left( v\left( \cdot \right) \right) ,
\end{equation*}%
where $J\left( \cdot \right) $ is its cost function as follows: 
\begin{equation}
J\left( v\left( \cdot \right) \right) =\mathbb{E}\left[ \int_{0}^{T}l\left(
s,x,u\left( s,x\right) ,\left( \nabla u\sigma \right) \left( s,x,u\right)
,v\left( s\right) \right) \text{d}s+\gamma \left( u\left( 0,x\right) \right) %
\right] .  \tag{1.4}
\end{equation}

As we have known, stochastic control problem of the SPDEs arising from
partial observation control has been studied by Mortensen [10], using a
dynamic programming approach, and subsequently by Bensoussan [2], [3], using
a maximum principle method. See [4], [16] and the references therein for
more information. Our approach differs from the one of Bensoussan. More
precisely, we relate the FBDSDE to one kind of SPDEs with control variables
where the control systems of SPDEs can be transformed to the relevant
control systems of FBDSDE. To our knowledge, this is the first time to treat
the optimal control problems of SPDEs from a new perspective of FBDSDE. It
is worth mentioning that the quasilinear SPDEs in [13] Øksendal considered
can just be related to our partially coupled FBDSDE. Recently, Zhang and Shi
[26], obtained the similar results, however, in their paper, the
coefficients $\sigma $ and $g$ do not contain the control variable,
respectively. The similar result for BDSDEs can be seen in [7].

This paper is organized as follows. Section 2 is devoting to stating the
problems and some assumptions. In Section 3 and Section 4, we give the
necessary and sufficient maximum principles for fully couple
forward-backward doubly stochastic control systems, respectively, in global
form. As an application, we study the optimal control of SPDEs in Section 5.
Finally, in Section 6 our results are further illustrated by solving optimal
controls of LQ problem and a special SPDEs using the Malliavin calculus,
respectively.

\section{Statement of the problems}

Let $\left( \Omega ,\mathcal{F},P\right) $ be a completed probability space, 
$\left\{ W_{t}\right\} _{t\geq 0}$ and $\left\{ B_{t}\right\} _{t\geq 0}$ be
two mutually independent standard Brownian motions, with value respectively
in $\mathbb{R}^{d}$ and $\mathbb{R}^{l},$ defined on $\left( \Omega ,%
\mathcal{F},P\right) $. Let $\mathcal{N}$ denote the class of $P$-null sets
of $\mathcal{F}$. For each $t\in \left[ 0,T\right] ,$ we define 
\begin{equation*}
\mathcal{F}_{t}^{W}\doteq \sigma \left\{ W_{r};\text{ }0\leq r\leq t\right\}
\bigvee \mathcal{N},\text{ \quad }\mathcal{F}_{t,T}^{B}\doteq \sigma \left\{
B_{r}-B_{t};\text{ }t\leq r\leq T\right\} \bigvee \mathcal{N},
\end{equation*}%
and 
\begin{equation*}
\mathcal{F}_{t}\doteq \mathcal{F}_{t}^{W}\bigvee \mathcal{F}_{t,T}^{B},\text{
}\forall t\in \left[ 0,T\right] .
\end{equation*}%
Note that $\left\{ \mathcal{F}_{t}^{W};t\in \left[ 0,T\right] \right\} $ is
an increasing filtration and $\left\{ \mathcal{F}_{t,T}^{B};t\in \left[ 0,T%
\right] \right\} $ is a decreasing filtration, and the collection $\left\{ 
\mathcal{F}_{t},t\in \left[ 0,T\right] \right\} $ is neither increasing nor
decreasing.

We denote $M^{2}\left( 0,T;\mathbb{R}^{n}\right) $ the space of (class of $%
dP\otimes dt$ a.e equal) all $\left\{ \mathcal{F}_{t}\right\} $-measurable $n
$-dimensional processes $\upsilon $ with norm of $\parallel \upsilon
\parallel _{M}\doteq \left[ \mathbb{E}\int_{0}^{T}|\upsilon (s)|^{2}\text{d}s%
\right] ^{\frac{1}{2}}<\infty .$ Obviously $M^{2}\left( 0,T;\mathbb{R}%
^{n}\right) $ is a Hilbert space. For any given $u\in M^{2}\left( 0,T\text{ }%
;\mathbb{R}^{n}\right) $ and $\upsilon \in M^{2}\left( 0,T\text{ };\mathbb{R}%
^{n}\right) ,$ one can define the (standard) forward It\^{o}'s integral $%
\int_{0}^{\cdot }u_{s}$d$\overrightarrow{W_{s}}$ and backward It\^{o}'s
integral $\int_{\cdot }^{T}\upsilon _{s}$d$\overleftarrow{B_{s}}$. They are
both in $M^{2}\left( 0,T;\mathbb{R}^{n}\right) ,$ (see [15] for details).

Let $L^{2}\left( \Omega ,\mathcal{F}_{T},P;\mathbb{R}^{n}\right) $ denote
the space of all $\left\{ \mathcal{F}_{T}\right\} $-measurable $\mathbb{R}%
^{n}$-valued random variable $\xi $ satisfying $\mathbb{E}\left\vert \xi
\right\vert ^{2}<\infty .$

\begin{definition}
A stochastic process $X=\left\{ X_t;t\geq 0\right\} $ is called $\mathcal{F}%
_t$-progressively measurable, if for any $t\geq 0$, $X$ on $\Omega \times %
\left[ 0,t\right] $ is measurable with respect to $\left( \mathcal{F}%
_t^W\times \mathcal{B}\left( \left[ 0,t\right] \right) \right) \vee \left( 
\mathcal{F}_{t,T}^B\times \mathcal{B}\left( \left[ t,T\right] \right)
\right) $.
\end{definition}

Under this framework, we consider the following forward-backward doubly
stochastic control system%
\begin{equation}
\left\{ 
\begin{array}{l}
\text{d}y\left( t\right) =f\left( t,y\left( t\right) ,Y\left( t\right)
,z\left( t\right) ,Z\left( t\right) ,v\left( t\right) \right) \text{d}t \\ 
\qquad +g\left( t,y\left( t\right) ,Y\left( t\right) ,z\left( t\right)
,Z\left( t\right) ,v\left( t\right) \right) \text{d}\overrightarrow{W_{t}}%
-z\left( t\right) \text{d}\overleftarrow{B_{t}}, \\ 
\text{d}Y\left( t\right) =-F\left( t,y\left( t\right) ,Y\left( t\right)
,z\left( t\right) ,Z\left( t\right) ,v\left( t\right) \right) \text{d}t \\ 
\qquad -G\left( t,y\left( t\right) ,Y\left( t\right) ,z\left( t\right)
,Z\left( t\right) ,v\left( t\right) \right) \text{d}\overleftarrow{B_{t}}%
+Z\left( t\right) \text{d}\overrightarrow{W_{t}}, \\ 
y\left( 0\right) =x_{0},\quad Y\left( T\right) =\varphi \left( y\left(
T\right) \right) ,%
\end{array}%
\right.   \tag{2.1}
\end{equation}%
where $\left( y\left( t\right) ,Y\left( t\right) ,z\left( t\right) ,Z\left(
t\right) ,v\left( t\right) \right) \in \mathbb{R}^{n}\mathbb{\times R}^{n}%
\mathbb{\times R}^{n\times l}\mathbb{\times R}^{n\times d}\mathbb{\times R}%
^{k},$ $x_{0}\in \mathbb{R}^{n}\mathbf{,}$ is a given constant$,$ $t>0$ and $%
T>0,$%
\begin{equation*}
\begin{array}{llll}
F: & \left[ 0,T\right] \times \mathbb{R}^{n}\mathbb{\times R}^{n}\mathbb{%
\times R}^{n\times l}\mathbb{\times R}^{n\times d}\mathbb{\times R}^{k} & 
\rightarrow  & \mathbb{R}^{n}\mathbf{,} \\ 
f: & \left[ 0,T\right] \times \mathbb{R}^{n}\mathbb{\times R}^{n}\mathbb{%
\times R}^{n\times l}\mathbb{\times R}^{n\times d}\mathbb{\times R}^{k} & 
\rightarrow  & \mathbb{R}^{n}\mathbf{,} \\ 
G: & \left[ 0,T\right] \times \mathbb{R}^{n}\mathbb{\times R}^{n}\mathbb{%
\times R}^{n\times l}\mathbb{\times R}^{n\times d}\mathbb{\times R}^{k} & 
\rightarrow  & \mathbb{R}^{n}\mathbf{,} \\ 
g: & \left[ 0,T\right] \times \mathbb{R}^{n}\mathbb{\times R}^{n}\mathbb{%
\times R}^{n\times l}\mathbb{\times R}^{n\times d}\mathbb{\times R}^{k} & 
\rightarrow  & \mathbb{R}^{n}\mathbf{,} \\ 
\varphi : & \mathbb{R}^{n}\rightarrow \mathbb{R}^{n}\mathbf{.} &  & 
\end{array}%
\end{equation*}%
Let $\mathcal{U}$ be a nonempty convex subset of $\mathbb{R}^{k}\mathbf{.}$
We define the admissible control set 
\begin{equation*}
\mathcal{U}_{ad}\doteq \left\{ v\left( \cdot \right) \in M^{2}\left( 0,T;%
\mathbb{R}^{k}\right) ;\text{ }v\left( t\right) \in \mathcal{U},\text{ }%
0\leq t\leq T,\text{ a.e., a.s.}\right\} .
\end{equation*}%
Our optimal control problem is to minimize the cost function: 
\begin{equation}
J\left( v\left( \cdot \right) \right) \doteq \mathbb{E}\left[
\int_{0}^{T}l\left( t,y\left( t\right) ,Y\left( t\right) ,z\left( t\right)
,Z\left( t\right) ,v\left( t\right) \right) \text{d}t+\Phi \left( y\left(
T\right) \right) +\gamma \left( Y\left( 0\right) \right) \right]   \tag{2.2}
\end{equation}%
over $\mathcal{U}_{ad}$, where 
\begin{eqnarray*}
l &:&\left[ 0,T\right] \times \mathbb{R}^{n}\mathbb{\times R}^{n}\mathbb{%
\times R}^{n\times l}\mathbb{\times R}^{n\times d}\mathbb{\times R}%
^{k}\rightarrow \mathbb{R}\mathbf{,} \\
\Phi  &:&\mathbb{R}^{n}\mathbf{\rightarrow }\mathbb{R}\mathbf{,} \\
\gamma  &:&\mathbb{R}^{n}\mathbf{\rightarrow }\mathbb{R}\mathbf{.}
\end{eqnarray*}%
An admissible control $u\left( \cdot \right) $ is called an optimal control
if it attains the minimum over $\mathcal{U}_{ad}$. That is to say, we want
to find a $u\left( \cdot \right) ,$ such that 
\begin{equation*}
J\left( u\left( \cdot \right) \right) \doteq \underset{v\left( \cdot \right)
\in \mathcal{U}_{ad}}{\inf }J\left( v\left( \cdot \right) \right) .
\end{equation*}%
(2.1) is called the state equation, the solution $\left( y\left( \cdot
\right) ,Y\left( \cdot \right) ,z\left( \cdot \right) ,Z\left( \cdot \right)
\right) $ corresponding to $u\left( \cdot \right) $ is called the optimal
trajectory. Next we will give some notations: 
\begin{equation*}
\zeta =\left( 
\begin{array}{c}
y \\ 
Y \\ 
z \\ 
Z%
\end{array}%
\right) ,\quad A\left( t,\zeta \right) =\left( 
\begin{array}{c}
-F \\ 
f \\ 
-G \\ 
g%
\end{array}%
\right) \left( t,\zeta \right) .
\end{equation*}%
We use the usual inner product $\left\langle \cdot ,\cdot \right\rangle $
and Euclidean norm $\left\vert \cdot \right\vert $ in $\mathbb{R}^{n},$ $%
\mathbb{R}^{n\times l}$, and $\mathbb{R}^{n\times d}.$ All the equalities
and inequalities mentioned in this paper are in the sense of $dt\otimes dP$
almost surely on $\left[ 0,T\right] \times \Omega .$ We assume that

\begin{enumerate}
\item[\textbf{(H1)}] Assume that%
\begin{equation*}
\left\{ 
\begin{array}{l}
\text{For each }\zeta \in \mathbb{R}^{n+n+n\times l+n\times d},\text{ }%
A\left( \cdot ,\zeta \right) \text{ is an }\mathcal{F}_{t}\text{-measurable
process defined on }\left[ 0,T\right] \text{ } \\ 
\text{with }A\left( \cdot ,0\right) \in M^{2}\left( 0,T;\mathbb{R}%
^{n+n+n\times l+n\times d}\right) .%
\end{array}%
\right. 
\end{equation*}

\item[\textbf{(H2)}] $A\left( t,\zeta \right) $ and $\varphi \left( y\right) 
$ satisfy Lipschitz conditions: there exists a constant $k>0,$ such that 
\begin{equation*}
\left\{ 
\begin{array}{l}
\left\vert A\left( t,\zeta \right) -A\left( t,\bar{\zeta}\right) \right\vert
\leq k\left\vert \zeta -\bar{\zeta}\right\vert ,\text{\quad }\forall \zeta ,%
\text{ }\bar{\zeta}\in \mathbb{R}^{n+n+n\times l+n\times d},\text{ }\forall
t\in \left[ 0,T\right] , \\ 
\left\vert \varphi \left( y\right) -\varphi \left( \bar{y}\right)
\right\vert \leq k\left\vert y-\bar{y}\right\vert ,\quad \forall y,\text{ }%
\bar{y}\in \mathbb{R}^{n}\mathbf{.}%
\end{array}%
\right. 
\end{equation*}
\end{enumerate}

The following monotonic conditions introduced in [18], are the main
assumptions in this paper.

\begin{enumerate}
\item[\textbf{(H3)}] Assume that%
\begin{equation*}
\left\{ 
\begin{array}{l}
\left\langle A\left( t,\zeta \right) -A\left( t,\bar{\zeta}\right) ,\zeta -%
\bar{\zeta}\right\rangle \leq -\mu \left\vert \zeta -\bar{\zeta}\right\vert
^{2}, \\ 
\quad \forall \zeta =\left( y,Y,z,Z\right) ^{T},\text{ }\bar{\zeta}=\left( 
\bar{y},\bar{Y},\bar{z},\bar{Z}\right) ^{T}\in \mathbb{R}^{n}\mathbb{\times R%
}^{n}\mathbb{\times R}^{n\times l}\mathbb{\times R}^{n\times d},\text{ }%
\forall t\in \left[ 0,T\right] . \\ 
\left\langle \varphi \left( y\right) -\varphi \left( \bar{y}\right) ,y-\bar{y%
}\right\rangle \geq 0,\text{ }\forall y,\text{ }\bar{y}\in \mathbb{R}^{n}%
\mathbf{.}%
\end{array}%
\right. 
\end{equation*}
\end{enumerate}

\noindent or

\begin{enumerate}
\item[\textbf{(H'3)}] Assume that%
\begin{equation*}
\left\{ 
\begin{array}{l}
\left\langle A\left( t,\zeta \right) -A\left( t,\bar{\zeta}\right) ,\zeta -%
\bar{\zeta}\right\rangle \geq \mu \left\vert \zeta -\bar{\zeta}\right\vert
^{2}, \\ 
\quad \forall \zeta =\left( y,Y,z,Z\right) ^{T},\text{ }\bar{\zeta}=\left( 
\bar{y},\bar{Y},\bar{z},\bar{Z}\right) ^{T}\in \mathbb{R}^{n}\mathbb{\times R%
}^{n}\mathbb{\times R}^{n\times l}\mathbb{\times R}^{n\times d},\text{ }%
\forall t\in \left[ 0,T\right] . \\ 
\left\langle \varphi \left( y\right) -\varphi \left( \bar{y}\right) ,y-\bar{y%
}\right\rangle \leq 0,\text{ }\forall y,\text{ }\bar{y}\in \mathbb{R}^{n}%
\mathbf{,}%
\end{array}%
\right. 
\end{equation*}

\noindent where $\mu $ is some positive constant.
\end{enumerate}

\begin{proposition}
For any given admissible control $v\left( \cdot \right) ,$ we assume (H1),
(H2) and (H3) (or (H1), (H2) and (H3)') hold. Then FBDSDE (2.1) has the
unique solution 
\begin{equation*}
\left( y\left( \cdot \right) ,Y\left( \cdot \right) ,z\left( \cdot \right)
,Z\left( \cdot \right) \right) \in M^{2}\left( 0,T;\mathbb{R}^{n}\mathbb{%
\times R}^{n}\mathbb{\times R}^{n\times l}\mathbb{\times R}^{n\times
d}\right) .
\end{equation*}
\end{proposition}

\noindent The proof of Proposition 2 can be seen in [18]. We assume:

\begin{enumerate}
\item[\textbf{(H4)}] 
\begin{equation*}
\left\{ 
\begin{array}{l}
\text{i) }F,\text{ }f,\text{ }G,\text{ }g,\text{ }\varphi ,\text{ }l,\text{ }%
\Phi ,\text{ }\gamma \text{ are continuously differentiable} \\ 
\quad \text{with respect to }\left( y,Y,z,Z,v\right) ,\text{ }y\text{ and }Y;
\\ 
\text{ii) The derivatives of }F,\text{ }f,\text{ }G,\text{ }g,\text{ }%
\varphi \text{ are bounded;} \\ 
\text{iii) The derivatives of }l\text{ are bounded by }C\left( 1+\left\vert
y\right\vert +\left\vert Y\right\vert +\left\vert z\right\vert +\left\vert
Z\right\vert +\left\vert v\right\vert \right) ; \\ 
\text{iv) The derivatives of }\Phi \text{ and }\gamma \text{ with respect to 
}y,\text{ }Y\text{ are bounded by } \\ 
\quad C\left( 1+\left\vert y\right\vert \right) \text{ and }C\left(
1+\left\vert Y\right\vert \right) ,\text{ respectively.}%
\end{array}%
\right. 
\end{equation*}
\end{enumerate}

\noindent Lastly, we need the following extension of It\^{o}'s formula (for
more details see [15]).

\begin{proposition}
Let $\alpha \in S^{2}\left( [0,T];\mathbb{R}^{k}\right) ,$ $\beta \in
M^{2}\left( [0,T];\mathbb{R}^{k}\right) ,$ $\gamma \in M^{2}\left( [0,T];%
\mathbb{R}^{k\times l}\right) ,$ $\delta \in S^{2}\left( [0,T];\mathbb{R}%
^{k\times d}\right) $ satisfy: 
\begin{equation*}
\alpha _{t}=\alpha _{0}+\int_{0}^{t}\beta _{s}ds+\int_{0}^{t}\gamma _{s}d%
\overleftarrow{B_{s}}+\int_{0}^{t}\delta _{s}d\overrightarrow{W_{s}},\text{ }%
0\leq t\leq T.
\end{equation*}%
Then 
\begin{eqnarray*}
\left\vert \alpha _{t}\right\vert ^{2} &=&\left\vert \alpha _{0}\right\vert
^{2}+2\int_{0}^{t}\left\langle \alpha _{s},\beta _{s}\right\rangle \text{d}%
s+2\int_{0}^{t}\left\langle \alpha _{s},\gamma _{s}\text{d}\overleftarrow{%
B_{s}}\right\rangle +2\int_{0}^{t}\left\langle \alpha _{s},\delta _{s}\text{d%
}\overrightarrow{W_{s}}\right\rangle \\
&&-\int_{0}^{t}\left\vert \gamma _{s}\right\vert ^{2}\text{d}%
s+\int_{0}^{t}\left\vert \delta _{s}\right\vert ^{2}\text{d}s, \\
\mathbb{E}\left\vert \alpha _{t}\right\vert ^{2} &=&\mathbb{E}\left\vert
\alpha _{0}\right\vert ^{2}+2\mathbb{E}\int_{0}^{t}\left\langle \alpha
_{s},\beta _{s}\right\rangle \text{d}s-\mathbb{E}\int_{0}^{t}\left\vert
\gamma _{s}\right\vert ^{2}\text{d}s+\mathbb{E}\int_{0}^{t}\left\vert \delta
_{s}\right\vert ^{2}\text{d}s.
\end{eqnarray*}%
More generally, if $\phi \in C^{2}\left( \mathbb{R}^{k}\right) ,$%
\begin{eqnarray*}
\phi \left( \alpha _{t}\right) &=&\phi \left( \alpha _{0}\right)
+\int_{0}^{t}\left\langle \phi ^{^{\prime }}\left( \alpha _{s}\right) ,\beta
_{s}\right\rangle \text{d}s+\int_{0}^{t}\left\langle \phi ^{^{\prime
}}\left( \alpha _{s}\right) ,\gamma _{s}\text{d}\overleftarrow{B_{s}}%
\right\rangle +\int_{0}^{t}\left\langle \phi ^{^{\prime }}\left( \alpha
_{s}\right) ,\delta _{s}\text{d}\overrightarrow{W_{s}}\right\rangle \\
&&\ \ \ -\frac{1}{2}\int_{0}^{t}Tr\left[ \phi ^{^{\prime \prime }}\left(
\alpha _{s}\right) \gamma _{s}\gamma _{s}^{\ast }\right] \text{d}s+\frac{1}{2%
}\int_{0}^{t}Tr\left[ \phi ^{^{\prime \prime }}\left( \alpha _{s}\right)
\delta _{s}\delta _{s}^{\ast }\right] \text{d}s.
\end{eqnarray*}
\end{proposition}

\noindent Here $S^{2}\left( 0,T;\mathbb{R}^{k}\right) $ denotes the space of
(classes of $dt\otimes dP$ a.e. equal) all $\mathcal{F}_{t}$-progressively
measurable $k$-dimensional processes $v$ with 
\begin{equation*}
\mathbb{E}\left( \sup\limits_{0\leq t\leq T}\left\vert v(t)\right\vert
^{2}\right) <\infty .
\end{equation*}

\section{A necessary maximum principle for optimal controls of
forward-backward doubly stochastic control systems}

We consider the forward-backward doubly stochastic control system (2.1) and
the cost function (2.2). Let $u\left( \cdot \right) $ be an optimal control
and $\left( y\left( \cdot \right) ,Y\left( \cdot \right) ,z\left( \cdot
\right) ,Z\left( \cdot \right) \right) $ be the corresponding trajectory.
Let $v\left( \cdot \right) $ be any given admissible control such that $%
u\left( \cdot \right) +v\left( \cdot \right) \in \mathcal{U}_{ad}.$ Since $%
\mathcal{U}_{ad}$ is convex, then for any $0\leq \rho \leq 1,$ $u_{\rho
}\left( \cdot \right) =u\left( \cdot \right) +\rho v\left( \cdot \right) $
is also in $\mathcal{U}_{ad}.$ The following technique is mainly borrowed
from [22].

We introduce the following variational equation of FBDSDE (2.1): 
\begin{equation}
\left\{ 
\begin{array}{l}
\text{d}y^1\left( t\right) =[f_y\left( t,y\left( t\right) ,Y\left( t\right)
,z\left( t\right) ,Z\left( t\right) ,u\left( t\right) \right) y^1\left(
t\right) \\ 
\qquad \quad +f_Y\left( t,y\left( t\right) ,Y\left( t\right) ,z\left(
t\right) ,Z\left( t\right) ,u\left( t\right) \right) Y^1\left( t\right) \\ 
\qquad \quad +f_z\left( t,y\left( t\right) ,Y\left( t\right) ,z\left(
t\right) ,Z\left( t\right) ,u\left( t\right) \right) z^1\left( t\right) \\ 
\qquad \quad +f_Z\left( t,y\left( t\right) ,Y\left( t\right) ,z\left(
t\right) ,Z\left( t\right) ,u\left( t\right) \right) Z^1\left( t\right) \\ 
\qquad \quad +f_v\left( t,y\left( t\right) ,Y\left( t\right) ,z\left(
t\right) ,Z\left( t\right) ,u\left( t\right) \right) v\left( t\right) ]\text{%
d}t \\ 
\qquad \quad +[g_y\left( t,y\left( t\right) ,Y\left( t\right) ,z\left(
t\right) ,Z\left( t\right) ,u\left( t\right) \right) y^1\left( t\right) \\ 
\qquad \quad +g_Y\left( t,y\left( t\right) ,Y\left( t\right) ,z\left(
t\right) ,Z\left( t\right) ,u\left( t\right) \right) Y^1\left( t\right) \\ 
\qquad \quad +g_z\left( t,y\left( t\right) ,Y\left( t\right) ,z\left(
t\right) ,Z\left( t\right) ,u\left( t\right) \right) z^1\left( t\right) \\ 
\qquad \quad +g_Z\left( t,y\left( t\right) ,Y\left( t\right) ,z\left(
t\right) ,Z\left( t\right) ,u\left( t\right) \right) Z^1\left( t\right) \\ 
\qquad \quad +g_v\left( t,y\left( t\right) ,Y\left( t\right) ,z\left(
t\right) ,Z\left( t\right) ,u\left( t\right) \right) v\left( t\right) ]\text{%
d}\overrightarrow{W_t}-z^1\left( t\right) \text{d}\overleftarrow{B_t}, \\ 
\text{d}Y^1\left( t\right) =-[F_y\left( t,y\left( t\right) ,Y\left( t\right)
,z\left( t\right) ,Z\left( t\right) ,u\left( t\right) \right) y^1\left(
t\right) \\ 
\qquad \quad +F_Y\left( t,y\left( t\right) ,Y\left( t\right) ,z\left(
t\right) ,Z\left( t\right) ,u\left( t\right) \right) Y^1\left( t\right) \\ 
\qquad \quad +F_z\left( t,y\left( t\right) ,Y\left( t\right) ,z\left(
t\right) ,Z\left( t\right) ,u\left( t\right) \right) z^1\left( t\right) \\ 
\qquad \quad +F_Z\left( t,y\left( t\right) ,Y\left( t\right) ,z\left(
t\right) ,Z\left( t\right) ,u\left( t\right) \right) Z^1\left( t\right) \\ 
\qquad \quad +F_v\left( t,y\left( t\right) ,Y\left( t\right) ,z\left(
t\right) ,Z\left( t\right) ,u\left( t\right) \right) v\left( t\right) ]\text{%
d}t \\ 
\qquad \quad -[G_y\left( t,y\left( t\right) ,Y\left( t\right) ,z\left(
t\right) ,Z\left( t\right) ,u\left( t\right) \right) y^1\left( t\right) \\ 
\qquad \quad +G_Y\left( t,y\left( t\right) ,Y\left( t\right) ,z\left(
t\right) ,Z\left( t\right) ,u\left( t\right) \right) Y^1\left( t\right) \\ 
\qquad \quad +G_z\left( t,y\left( t\right) ,Y\left( t\right) ,z\left(
t\right) ,Z\left( t\right) ,u\left( t\right) \right) z^1\left( t\right) \\ 
\qquad \quad +G_Z\left( t,y\left( t\right) ,Y\left( t\right) ,z\left(
t\right) ,Z\left( t\right) ,u\left( t\right) \right) Z^1\left( t\right) \\ 
\qquad \quad +G_v\left( t,y\left( t\right) ,Y\left( t\right) ,z\left(
t\right) ,Z\left( t\right) ,u\left( t\right) \right) v\left( t\right) ]\text{%
d}\overleftarrow{B_t}+Z^1\left( t\right) \text{d}\overrightarrow{W_t}, \\ 
y^1\left( 0\right) =0,\quad Y^1\left( t\right) =\varphi _y\left( y\left(
T\right) \right) y^1\left( T\right) .%
\end{array}
\right.  \tag{3.1}
\end{equation}
From (H3), (H4) and Proposition 2, it is easy to check that (3.1) satisfies
(H1), (H2) and (H3). Then there exists a unique quadruple of $\left(
y^1\left( t\right) ,Y^1\left( t\right) ,z^1\left( t\right) ,Z^1\left(
t\right) \right) $ in $M^2\left( 0,T\right) $ satisfying FBDSDE (3.1). We
denote by $\left( y_\rho \left( t\right) ,Y_\rho \left( t\right) ,z_\rho
\left( t\right) ,Z_\rho \left( t\right) \right) $ the trajectory of FBDSDE
(2.1) corresponding to $u_\rho \left( \cdot \right) $ as follows$.$%
\begin{equation*}
\left\{ 
\begin{array}{l}
\text{d}y_\rho \left( t\right) =f\left( t,y_\rho \left( t\right) ,Y_\rho
\left( t\right) ,z_\rho \left( t\right) ,Z_\rho \left( t\right) ,u_\rho
\left( t\right) \right) \text{d}t \\ 
\qquad +g\left( t,y_\rho \left( t\right) ,Y_\rho \left( t\right) ,z_\rho
\left( t\right) ,Z_\rho \left( t\right) ,u_\rho \left( t\right) \right) 
\text{d}\overrightarrow{W_t}-z_\rho \left( t\right) \text{d}\overleftarrow{%
B_t}, \\ 
\text{d}Y_\rho \left( t\right) =-F\left( t,y_\rho \left( t\right) ,Y_\rho
\left( t\right) ,z_\rho \left( t\right) ,Z_\rho \left( t\right) ,u_\rho
\left( t\right) \right) \text{d}t \\ 
\qquad -G\left( t,y_\rho \left( t\right) ,Y_\rho \left( t\right) ,z_\rho
\left( t\right) ,Z_\rho \left( t\right) ,u_\rho \left( t\right) \right) 
\text{d}\overleftarrow{B_t}+Z_\rho \left( t\right) \text{d}\overrightarrow{%
W_t}, \\ 
y_\rho \left( 0\right) =x_0,\quad Y_\rho \left( T\right) =\varphi \left(
y_\rho \left( T\right) \right) ,%
\end{array}
\right.
\end{equation*}
Then we will study the solutions to forward-backward doubly stochastic
control systems with parameter.

\begin{lemma}
Assume that (H1)-(H4) hold. Then we have 
\begin{eqnarray*}
\underset{\rho \rightarrow 0}{\lim }\frac{y_\rho \left( t\right) -y\left(
t\right) }\rho =y^1\left( t\right) , \\
\underset{\rho \rightarrow 0}{\lim }\frac{Y_\rho \left( t\right) -Y\left(
t\right) }\rho =Y^1\left( t\right) , \\
\underset{\rho \rightarrow 0}{\lim }\frac{z_\rho \left( t\right) -z\left(
t\right) }\rho =z^1\left( t\right) , \\
\underset{\rho \rightarrow 0}{\lim }\frac{Z_\rho \left( t\right) -Z\left(
t\right) }\rho =Z^1\left( t\right) ,
\end{eqnarray*}
where the limits are in $M^2\left( 0,T\right) $.
\end{lemma}

\begin{proof}
Firstly, we show the continuous dependence of solutions with respect to the
parameter $\rho $. Let
\begin{eqnarray*}
\hat{y}\left( t\right) &=&y_{\rho }\left( t\right) -y\left( t\right) , \\
\hat{Y}\left( t\right) &=&Y_{\rho }\left( t\right) -Y\left( t\right) , \\
\hat{z}\left( t\right) &=&z_{\rho }\left( t\right) -z\left( t\right) , \\
\hat{Z}\left( t\right) &=&Z_{\rho }\left( t\right) -Z\left( t\right) .
\end{eqnarray*}%
We have
\begin{equation*}
\left\{
\begin{array}{l}
\text{d}\hat{y}\left( t\right) =[f(t,y_{\rho }\left( t\right) ,Y_{\rho
}\left( t\right) ,z_{\rho }\left( t\right) ,Z_{\rho }\left( t\right)
,u\left( t\right) +\rho v\left( t\right) ) \\
\quad \qquad -f(t,y_{\rho }\left( t\right) ,Y_{\rho }\left( t\right)
,z_{\rho }\left( t\right) ,Z_{\rho }\left( t\right) ,u\left( t\right) ) \\
\quad \qquad +f(t,y_{\rho }\left( t\right) ,Y_{\rho }\left( t\right)
,z_{\rho }\left( t\right) ,Z_{\rho }\left( t\right) ,u\left( t\right) ) \\
\quad \qquad -f(t,y\left( t\right) ,Y\left( t\right) ,z\left( t\right)
,Z\left( t\right) ,u\left( t\right) )]\text{d}t \\
\quad \qquad +[g(t,y_{\rho }\left( t\right) ,Y_{\rho }\left( t\right)
,z_{\rho }\left( t\right) ,Z_{\rho }\left( t\right) ,u\left( t\right) +\rho
v\left( t\right) ) \\
\quad \qquad -g(t,y_{\rho }\left( t\right) ,Y_{\rho }\left( t\right)
,z_{\rho }\left( t\right) ,Z_{\rho }\left( t\right) ,u\left( t\right) ) \\
\quad \qquad +g(t,y_{\rho }\left( t\right) ,Y_{\rho }\left( t\right)
,z_{\rho }\left( t\right) ,Z_{\rho }\left( t\right) ,u\left( t\right) ) \\
\quad \qquad -g(t,y\left( t\right) ,Y\left( t\right) ,z\left( t\right)
,Z\left( t\right) ,u\left( t\right) )]\text{d}\overrightarrow{W_{t}}-\hat{z}%
\left( t\right) \text{d}\overleftarrow{B_{t}}, \\
\text{d}\hat{Y}\left( t\right) =-[F(t,y_{\rho }\left( t\right) ,Y_{\rho
}\left( t\right) ,z_{\rho }\left( t\right) ,Z_{\rho }\left( t\right)
,u\left( t\right) +\rho v\left( t\right) ) \\
\quad \qquad -F(t,y_{\rho }\left( t\right) ,Y_{\rho }\left( t\right)
,z_{\rho }\left( t\right) ,Z_{\rho }\left( t\right) ,u\left( t\right) ) \\
\quad \qquad +F\left( t,y_{\rho }\left( t\right) ,Y_{\rho }\left( t\right)
,z_{\rho }\left( t\right) ,Z_{\rho }\left( t\right) ,u\left( t\right) \right)
\\
\quad \qquad -F(t,y\left( t\right) ,Y\left( t\right) ,z\left( t\right)
,Z\left( t\right) ,u\left( t\right) )]\text{d}t \\
\quad \qquad -[G(t,y_{\rho }\left( t\right) ,Y_{\rho }\left( t\right)
,z_{\rho }\left( t\right) ,Z_{\rho }\left( t\right) ,u\left( t\right) +\rho
v\left( t\right) ) \\
\quad \qquad -G(t,y_{\rho }\left( t\right) ,Y_{\rho }\left( t\right)
,z_{\rho }\left( t\right) ,Z_{\rho }\left( t\right) ,u\left( t\right) ) \\
\quad \qquad +G(t,y_{\rho }\left( t\right) ,Y_{\rho }\left( t\right)
,z_{\rho }\left( t\right) ,Z_{\rho }\left( t\right) ,u\left( t\right) ) \\
\quad \qquad -G(t,y\left( t\right) ,Y\left( t\right) ,z\left( t\right)
,Z\left( t\right) ,u\left( t\right) )]\text{d}\overleftarrow{B_{t}}+\hat{Z}%
\left( t\right) \text{d}\overrightarrow{W_{t}}, \\
\hat{y}\left( 0\right) =0,\quad \hat{Y}\left( T\right) =\varphi \left(
y_{\rho }\left( T\right) \right) -\varphi \left( y\left( T\right) \right) .%
\end{array}%
\right.
\end{equation*}%
We will prove $\left( \hat{y}\left( t\right) ,\hat{Y}\left( t\right) ,\hat{z}%
\left( t\right) ,\hat{Z}\left( t\right) \right) $ converge to $0$ in $%
M^{2}\left( 0,T\right) $ as $\rho \rightarrow 0.$ Applying Itô's formula to $%
\left\langle \hat{y}\left( t\right) ,\hat{Y}\left( t\right) \right\rangle $
on $\left[ 0,T\right] ,$ and by (H4) it follows that
\begin{eqnarray*}
&&\mathbb{E}\left\langle \hat{y}\left( T\right) ,\varphi \left( y_{\rho
}\left( T\right) \right) -\varphi \left( y\left( T\right) \right)
\right\rangle \\
&=&\mathbb{E}\int_{0}^{T}\left\langle A\left( t,\xi _{\rho }\right) -A\left(
t,\xi \right) ,\xi _{\rho }-\xi \right\rangle \text{d}t \\
&&-\mathbb{E}\int_{0}^{T}\hat{y}\left( t\right) [F\left( t,y_{\rho }\left(
t\right) ,Y_{\rho }\left( t\right) ,z_{\rho }\left( t\right) ,Z_{\rho
}\left( t\right) ,u\left( t\right) +\rho v\left( t\right) \right) \\
&&-F\left( t,y_{\rho }\left( t\right) ,Y_{\rho }\left( t\right) ,z_{\rho
}\left( t\right) ,Z_{\rho }\left( t\right) ,u\left( t\right) \right) ]\text{d%
}t \\
&&+\mathbb{E}\int_{0}^{T}\hat{Y}\left( t\right) [f\left( t,y_{\rho }\left(
t\right) ,Y_{\rho }\left( t\right) ,z_{\rho }\left( t\right) ,Z_{\rho
}\left( t\right) ,u\left( t\right) +\rho v\left( t\right) \right) \\
&&-f\left( t,y_{\rho }\left( t\right) ,Y_{\rho }\left( t\right) ,z_{\rho
}\left( t\right) ,Z_{\rho }\left( t\right) ,u\left( t\right) \right) ]\text{d%
}t \\
&&-\mathbb{E}\int_{0}^{T}\hat{z}\left( t\right) [G\left( t,y_{\rho }\left(
t\right) ,Y_{\rho }\left( t\right) ,z_{\rho }\left( t\right) ,Z_{\rho
}\left( t\right) ,u\left( t\right) +\rho v\left( t\right) \right) \\
&&-G\left( t,y_{\rho }\left( t\right) ,Y_{\rho }\left( t\right) ,z_{\rho
}\left( t\right) ,Z_{\rho }\left( t\right) ,u\left( t\right) \right) ]\text{d%
}t \\
&&+\mathbb{E}\int_{0}^{T}\hat{Z}\left( t\right) [g\left( t,y_{\rho }\left(
t\right) ,Y_{\rho }\left( t\right) ,z_{\rho }\left( t\right) ,Z_{\rho
}\left( t\right) ,u\left( t\right) +\rho v\left( t\right) \right) \\
&&-g\left( t,y_{\rho }\left( t\right) ,Y_{\rho }\left( t\right) ,z_{\rho
}\left( t\right) ,Z_{\rho }\left( t\right) ,u\left( t\right) \right) ]\text{d%
}t \\
&\leq &-\mu \mathbb{E}\int_{0}^{T}\left[ \left\vert \hat{y}\left( t\right)
\right\vert ^{2}+\left\vert \hat{Y}\left( t\right) \right\vert
^{2}+\left\vert \hat{z}\left( t\right) \right\vert ^{2}+\left\vert \hat{Z}%
\left( t\right) \right\vert ^{2}\right] \text{d}t \\
&&+\frac{\mu }{4}\mathbb{E}\int_{0}^{T}\left[ \left\vert \hat{y}\left(
t\right) \right\vert ^{2}+\left\vert \hat{Y}\left( t\right) \right\vert
^{2}+\left\vert \hat{z}\left( t\right) \right\vert ^{2}+\left\vert \hat{Z}%
\left( t\right) \right\vert ^{2}\right] \text{d}t \\
&&+\frac{1}{\mu }\rho ^{2}C\mathbb{E}\int_{0}^{T}\left\vert v\left( t\right)
\right\vert ^{2}\text{d}t,
\end{eqnarray*}%
where
\begin{eqnarray*}
\xi _{\rho }\left( t\right) &=&\left( y_{\rho }\left( t\right) ,Y_{\rho
}\left( t\right) ,z_{\rho }\left( t\right) ,Z_{\rho }\left( t\right)
,u\left( t\right) \right) ^{T}, \\
\xi \left( t\right) &=&\left( y\left( t\right) ,Y\left( t\right) ,z\left(
t\right) ,Z\left( t\right) ,u\left( t\right) \right) ^{T}, \\
A\left( t,\xi \right) &=&\left(
\begin{array}{c}
-F\left( t,\xi \right) \\
f\left( t,\xi \right) \\
-G\left( t,\xi \right) \\
g\left( t,\xi \right)%
\end{array}%
\right) ,\text{ \quad }A\left( t,\xi _{\rho }\right) =\left(
\begin{array}{c}
-F\left( t,\xi _{\rho }\right) \\
f\left( t,\xi _{\rho }\right) \\
-G\left( t,\xi _{\rho }\right) \\
g\left( t,\xi _{\rho }\right)%
\end{array}%
\right) .
\end{eqnarray*}%
Thus we get
\begin{equation*}
\mathbb{E}\int_{0}^{T}\left[ \left\vert \hat{y}\left( t\right) \right\vert
^{2}+\left\vert \hat{Y}\left( t\right) \right\vert ^{2}+\left\vert \hat{z}%
\left( t\right) \right\vert ^{2}+\left\vert \hat{Z}\left( t\right)
\right\vert ^{2}\right] \text{d}t\leq \rho ^{2}C\mathbf{E}%
\int_{0}^{T}\left\vert v\left( t\right) \right\vert ^{2}\text{d}t.
\end{equation*}%
Then it follows that $\left( \hat{y}\left( t\right) ,\hat{Y}\left( t\right) ,%
\hat{z}\left( t\right) ,\hat{Z}\left( t\right) \right) $ converge to $0$ in $%
M^{2}\left( 0,T\right) $ as $\rho $ tends to $0.$ Set
\begin{eqnarray*}
\triangle y\left( t\right) &=&\frac{y_{\rho }\left( t\right) -y\left(
t\right) }{\rho }, \\
\triangle Y\left( t\right) &=&\frac{Y_{\rho }\left( t\right) -Y\left(
t\right) }{\rho }, \\
\triangle z\left( t\right) &=&\frac{z_{\rho }\left( t\right) -z\left(
t\right) }{\rho }, \\
\triangle Z\left( t\right) &=&\frac{Z_{\rho }\left( t\right) -Z\left(
t\right) }{\rho },
\end{eqnarray*}%
then
\begin{equation*}
\left\{
\begin{array}{l}
\text{d}\triangle y\left( t\right) =\frac{f\left( t,y_{\rho }\left( t\right)
,Y_{\rho }\left( t\right) ,z_{\rho }\left( t\right) ,Z_{\rho }\left(
t\right) ,u\left( t\right) +\rho v\left( t\right) \right) -f\left( t,y\left(
t\right) ,Y\left( t\right) ,z\left( t\right) ,Z\left( t\right) ,u\left(
t\right) \right) }{\rho }\text{d}t \\
\qquad \quad +\frac{g\left( t,y_{\rho }\left( t\right) ,Y_{\rho }\left(
t\right) ,z_{\rho }\left( t\right) ,Z_{\rho }\left( t\right) ,u\left(
t\right) +\rho v\left( t\right) \right) -g\left( t,y\left( t\right) ,Y\left(
t\right) ,z\left( t\right) ,Z\left( t\right) ,u\left( t\right) \right) }{%
\rho }\text{d}\overrightarrow{W_{t}} \\
\qquad \quad -\triangle z\left( t\right) \text{d}\overleftarrow{B_{t}}, \\
-\text{d}\triangle Y\left( t\right) =\frac{F\left( t,y_{\rho }\left(
t\right) ,Y_{\rho }\left( t\right) ,z_{\rho }\left( t\right) ,Z_{\rho
}\left( t\right) ,u\left( t\right) +\rho v\left( t\right) \right) -F\left(
t,y\left( t\right) ,Y\left( t\right) ,z\left( t\right) ,Z\left( t\right)
,u\left( t\right) \right) }{\rho }\text{d}t \\
\qquad \quad +\frac{G\left( t,y_{\rho }\left( t\right) ,Y_{\rho }\left(
t\right) ,z_{\rho }\left( t\right) ,Z_{\rho }\left( t\right) ,u\left(
t\right) +\rho v\left( t\right) \right) -G\left( t,y\left( t\right) ,Y\left(
t\right) ,z\left( t\right) ,Z\left( t\right) ,u\left( t\right) \right) }{%
\rho }\text{d}\overleftarrow{B_{t}} \\
\qquad \quad -\triangle Z\left( t\right) \text{d}\overrightarrow{W_{t}}, \\
\triangle y\left( 0\right) =0,\quad \triangle Y\left( T\right) =\frac{%
\varphi \left( y_{\rho }\left( T\right) -\varphi \left( y\left( T\right)
\right) \right) }{\rho }.%
\end{array}%
\right.
\end{equation*}%
The above equations can be expressed as follows
\begin{equation*}
\left\{
\begin{array}{l}
\text{d}\triangle y\left( t\right) =\bar{f}\left( t,\triangle y\left(
t\right) ,\triangle Y\left( t\right) ,\triangle z\left( t\right) ,\triangle
Z\left( t\right) ,v\left( t\right) \right) \text{d}t \\
\quad \qquad +\bar{g}\left( t,\triangle y\left( t\right) ,\triangle Y\left(
t\right) ,\triangle z\left( t\right) ,\triangle Z\left( t\right) ,v\left(
t\right) \right) \text{d}\overrightarrow{W_{t}} \\
\quad \qquad -\triangle z\left( t\right) \text{d}\overleftarrow{B_{t}}, \\
-\text{d}\triangle Y\left( t\right) =\bar{F}\left( t,\triangle y\left(
t\right) ,\triangle Y\left( t\right) ,\triangle z\left( t\right) ,\triangle
Z\left( t\right) ,v\left( t\right) \right) \text{d}t \\
\quad \qquad +\bar{G}\left( t,\triangle y\left( t\right) ,\triangle Y\left(
t\right) ,\triangle z\left( t\right) ,\triangle Z\left( t\right) ,v\left(
t\right) \right) \text{d}\overleftarrow{B_{t}} \\
\quad \qquad -\triangle Z\left( t\right) \text{d}\overrightarrow{W_{t}}, \\
\triangle y\left( 0\right) =0,\quad \triangle Y\left( T\right) =\frac{%
\varphi \left( y_{\rho }\left( T\right) \right) -\varphi \left( y\left(
T\right) \right) }{\rho },%
\end{array}%
\right.
\end{equation*}%
where $\bar{\theta}=\bar{f},$ $\bar{F}$, $\bar{g}$, $\bar{G}$, respectively,
\begin{equation*}
\bar{\theta}\left( t,\triangle y,\triangle Y,\triangle z,\triangle
Z,v\right) =A^{\theta }\left( t\right) \triangle y+B^{\theta }\left(
t\right) \triangle Y+C^{\theta }\left( t\right) \triangle z+D^{\theta
}\left( t\right) \triangle Z+E^{\theta }\left( t\right) v,
\end{equation*}%
and
\begin{eqnarray*}
A^{\theta }\left( t\right) &=&\left\{
\begin{array}{l}
\frac{\theta \left( t,y_{\rho }\left( t\right) ,Y_{\rho }\left( t\right)
,z_{\rho }\left( t\right) ,Z_{\rho }\left( t\right) ,u\left( t\right) +\rho
v\left( t\right) \right) -\theta \left( t,y\left( t\right) ,Y_{\rho }\left(
t\right) ,z_{\rho }\left( t\right) ,Z_{\rho }\left( t\right) ,u\left(
t\right) +\rho v\left( t\right) \right) }{y_{\rho }\left( t\right) -y\left(
t\right) },\quad y_{\rho }\left( t\right) -y\left( t\right) \neq 0, \\
0,\quad \quad \text{otherwise;}%
\end{array}%
\right. \\
B^{\theta }\left( t\right) &=&\left\{
\begin{array}{l}
\frac{\theta \left( t,y\left( t\right) ,Y_{\rho }\left( t\right) ,z_{\rho
}\left( t\right) ,Z_{\rho }\left( t\right) ,u\left( t\right) +\rho v\left(
t\right) \right) -\theta \left( t,y\left( t\right) ,Y\left( t\right)
,z_{\rho }\left( t\right) ,Z_{\rho }\left( t\right) ,u\left( t\right) +\rho
v\left( t\right) \right) }{Y_{\rho }\left( t\right) -Y\left( t\right) }%
,\quad Y_{\rho }\left( t\right) -Y\left( t\right) \neq 0, \\
0,\quad \quad \text{otherwise;}%
\end{array}%
\right. \\
C^{\theta }\left( t\right) &=&\left\{
\begin{array}{l}
\frac{\theta \left( t,y\left( t\right) ,Y\left( t\right) ,z_{\rho }\left(
t\right) ,Z_{\rho }\left( t\right) ,u\left( t\right) +\rho v\left( t\right)
\right) -\theta \left( t,y\left( t\right) ,Y\left( t\right) ,z\left(
t\right) ,Z_{\rho }\left( t\right) ,u\left( t\right) +\rho v\left( t\right)
\right) }{z_{\rho }\left( t\right) -z\left( t\right) },\quad z_{\rho }\left(
t\right) -z\left( t\right) \neq 0, \\
0,\quad \quad \text{otherwise;}%
\end{array}%
\right. \\
D^{\theta }\left( t\right) &=&\left\{
\begin{array}{l}
\frac{\theta \left( t,y\left( t\right) ,Y\left( t\right) ,z\left( t\right)
,Z_{\rho }\left( t\right) ,u\left( t\right) +\rho v\left( t\right) \right)
-\theta \left( t,y\left( t\right) ,Y\left( t\right) ,z\left( t\right)
,Z\left( t\right) ,u\left( t\right) +\rho v\left( t\right) \right) }{Z_{\rho
}\left( t\right) -Z\left( t\right) },\quad Z_{\rho }\left( t\right) -Z\left(
t\right) \neq 0, \\
0,\quad \quad \text{otherwise;}%
\end{array}%
\right. \\
E^{\theta }\left( t\right) &=&\left\{
\begin{array}{l}
\frac{\theta \left( t,y\left( t\right) ,Y\left( t\right) ,z\left( t\right)
,Z\left( t\right) ,u\left( t\right) +\rho v\left( t\right) \right) -\theta
\left( t,y\left( t\right) ,Y\left( t\right) ,z\left( t\right) ,Z\left(
t\right) ,u\left( t\right) \right) }{\rho v\left( t\right) },\quad \rho
v\left( t\right) \neq 0, \\
0,\quad \quad \text{otherwise.}%
\end{array}%
\right.
\end{eqnarray*}%
From the continuous dependence of solutions with respect to the parameter $%
\rho $, it follows that
\begin{eqnarray*}
\underset{\rho \rightarrow 0}{\lim }A^{\theta }\left( t\right) &=&\theta
_{y}\left( t,y\left( t\right) ,Y\left( t\right) ,z\left( t\right) ,Z\left(
t\right) ,u\left( t\right) \right) , \\
\underset{\rho \rightarrow 0}{\lim }B^{\theta }\left( t\right) &=&\theta
_{Y}\left( t,y\left( t\right) ,Y\left( t\right) ,z\left( t\right) ,Z\left(
t\right) ,u\left( t\right) \right) , \\
\underset{\rho \rightarrow 0}{\lim }C^{\theta }\left( t\right) &=&\theta
_{z}\left( t,y\left( t\right) ,Y\left( t\right) ,z\left( t\right) ,Z\left(
t\right) ,u\left( t\right) \right) , \\
\underset{\rho \rightarrow 0}{\lim }D^{\theta }\left( t\right) &=&\theta
_{Z}\left( t,y\left( t\right) ,Y\left( t\right) ,z\left( t\right) ,Z\left(
t\right) ,u\left( t\right) \right) , \\
\underset{\rho \rightarrow 0}{\lim }E^{\theta }\left( t\right) &=&\theta
_{v}\left( t,y\left( t\right) ,Y\left( t\right) ,z\left( t\right) ,Z\left(
t\right) ,u\left( t\right) \right) .
\end{eqnarray*}%
According to the continuous dependence of solutions with respect to the
parameter and the uniqueness of solutions of FBDSDE (3.1), the solutions $%
\left( \triangle y\left( t\right) ,\triangle Y\left( t\right) ,\triangle
z\left( t\right) ,\triangle Z\left( t\right) \right) $ converge to $\left(
y^{1}\left( t\right) ,Y^{1}\left( t\right) ,z^{1}\left( t\right)
,Z^{1}\left( t\right) \right) $ in $M^{2}\left( 0,T;\mathbb{R}^{n}\mathbb{%
\times R}^{n}\mathbb{\times R}^{n\times l}\mathbb{\times R}^{n\times
d}\right) $ as $\rho \rightarrow 0.$ The proof is completed.
\end{proof}

Now we give the variational inequality.

\begin{lemma}
Assume that (H1)-(H4) hold. Then we have 
\begin{eqnarray*}
&&\ \mathbb{E}\left[ \Phi _{y}\left( y\left( T\right) \right) y^{1}\left(
T\right) \right] +\mathbb{E}\left[ \gamma _{Y}\left( Y\left( 0\right)
\right) Y^{1}\left( 0\right) \right] \\
&&\ +\mathbb{E}\left[ \int_{0}^{T}[l_{y}\left( t,y\left( t\right) ,Y\left(
t\right) ,z\left( t\right) ,Z\left( t\right) ,u\left( t\right) \right)
y^{1}\left( t\right) \right. \\
&&\ +l_{Y}\left( t,y\left( t\right) ,Y\left( t\right) ,z\left( t\right)
,Z\left( t\right) ,u\left( t\right) \right) Y^{1}\left( t\right) \\
&&\ +l_{z}\left( t,y\left( t\right) ,Y\left( t\right) ,z\left( t\right)
,Z\left( t\right) ,u\left( t\right) \right) z^{1}\left( t\right) \\
&&\ +l_{Z}\left( t,y\left( t\right) ,Y\left( t\right) ,z\left( t\right)
,Z\left( t\right) ,u\left( t\right) \right) Z^{1}\left( t\right) \\
&&\ \left. +l_{v}\left( t,y\left( t\right) ,Y\left( t\right) ,z\left(
t\right) ,Z\left( t\right) ,u\left( t\right) \right) v\left( t\right) ]\text{%
d}t\right] \\
\ &\geq &0.
\end{eqnarray*}
\end{lemma}

\begin{proof}
From Lemma 4 and (H4), we can get
\begin{eqnarray*}
\underset{\rho \rightarrow 0}{\lim }\frac{\mathbb{E}\left[ \Phi \left(
y_{\rho }\left( T\right) \right) -\Phi \left( y\left( T\right) \right) %
\right] }{\rho } &=&\mathbb{E}\Phi _{y}\left( y\left( T\right) \right)
y^{1}\left( T\right) , \\
\underset{\rho \rightarrow 0}{\lim }\frac{\mathbb{E}\left[ \gamma \left(
Y_{\rho }\left( 0\right) \right) -\gamma \left( Y\left( 0\right) \right) %
\right] }{\rho } &=&\mathbb{E}\gamma _{Y}\left( Y\left( 0\right) \right)
Y^{1}\left( 0\right) ,
\end{eqnarray*}%
and
\begin{eqnarray*}
&&\underset{\rho \rightarrow 0}{\lim }\rho ^{-1}\mathbb{E}\left[
\int_{0}^{T}l\left( t,y_{\rho }\left( t\right) ,Y_{\rho }\left( t\right)
,z_{\rho }\left( t\right) ,Z_{\rho }\left( t\right) ,u\left( t\right) +\rho
v\left( t\right) \right) \right. \\
&&\left. -l\left( t,y\left( t\right) ,Y\left( t\right) ,z\left( t\right)
,Z\left( t\right) ,u\left( t\right) \right) \text{d}t\right] \\
&=&\mathbb{E}\int_{0}^{T}[l_{y}\left( t,y\left( t\right) ,Y\left( t\right)
,z\left( t\right) ,Z\left( t\right) ,u\left( t\right) \right) y^{1}\left(
t\right) \\
&&+l_{Y}\left( t,y\left( t\right) ,Y\left( t\right) ,z\left( t\right)
,Z\left( t\right) ,u\left( t\right) \right) Y^{1}\left( t\right) \\
&&+l_{z}\left( t,y\left( t\right) ,Y\left( t\right) ,z\left( t\right)
,Z\left( t\right) ,u\left( t\right) \right) z^{1}\left( t\right) \\
&&+l_{Z}\left( t,y\left( t\right) ,Y\left( t\right) ,z\left( t\right)
,Z\left( t\right) ,u\left( t\right) \right) Z^{1}\left( t\right) \\
&&+l_{v}\left( t,y\left( t\right) ,Y\left( t\right) ,z\left( t\right)
,Z\left( t\right) ,u\left( t\right) \right) v\left( t\right) ]\text{d}t.
\end{eqnarray*}%
On the other hand, since $u\left( \cdot \right) $ is an optimal control, it
follows that
\begin{equation*}
\rho ^{-1}\left[ J\left( u\left( \cdot \right) +\rho v\left( \cdot \right)
\right) -J\left( u\left( \cdot \right) \right) \right] \geq 0.
\end{equation*}%
Therefore the desired result is obtained.
\end{proof}

Now we introduce the adjoint equation by virtue of dual technique and
Hamilton function for our problem. From the variational inequality obtained
in Lemma 5, the maximum principle can be proved by using It\^{o}'s formula.
The adjoint equations are 
\begin{equation}
\left\{ 
\begin{array}{l}
\text{d}p\left( t\right) =[F_{Y}\left( t,y\left( t\right) ,Y\left( t\right)
,z\left( t\right) ,Z\left( t\right) ,u\left( t\right) \right) p\left(
t\right)  \\ 
\quad \qquad -f_{Y}\left( t,y\left( t\right) ,Y\left( t\right) ,z\left(
t\right) ,Z\left( t\right) ,u\left( t\right) \right) q\left( t\right)  \\ 
\quad \qquad +G_{Y}\left( t,y\left( t\right) ,Y\left( t\right) ,z\left(
t\right) ,Z\left( t\right) ,u\left( t\right) \right) k\left( t\right)  \\ 
\quad \qquad -g_{Y}\left( t,y\left( t\right) ,Y\left( t\right) ,z\left(
t\right) ,Z\left( t\right) ,u\left( t\right) \right) h\left( t\right)  \\ 
\quad \qquad -l_{Y}\left( t,y\left( t\right) ,Y\left( t\right) ,z\left(
t\right) ,Z\left( t\right) ,u\left( t\right) \right) ]\text{d}t \\ 
\quad \qquad +[F_{Z}\left( t,y\left( t\right) ,Y\left( t\right) ,z\left(
t\right) ,Z\left( t\right) ,u\left( t\right) \right) p\left( t\right)  \\ 
\quad \qquad -f_{Z}\left( t,y\left( t\right) ,Y\left( t\right) ,z\left(
t\right) ,Z\left( t\right) ,u\left( t\right) \right) q\left( t\right)  \\ 
\quad \qquad +G_{Z}\left( t,y\left( t\right) ,Y\left( t\right) ,z\left(
t\right) ,Z\left( t\right) ,u\left( t\right) \right) k\left( t\right)  \\ 
\quad \qquad -g_{Z}\left( t,y\left( t\right) ,Y\left( t\right) ,z\left(
t\right) ,Z\left( t\right) ,u\left( t\right) \right) h\left( t\right)  \\ 
\quad \qquad -l_{Z}\left( t,y\left( t\right) ,Y\left( t\right) ,z\left(
t\right) ,Z\left( t\right) ,u\left( t\right) \right) ]\text{d}%
\overrightarrow{W_{t}}-k_{t}\text{d}\overleftarrow{B_{t}}, \\ 
\text{d}q\left( t\right) =[F_{y}\left( t,y\left( t\right) ,Y\left( t\right)
,z\left( t\right) ,Z\left( t\right) ,u\left( t\right) \right) p\left(
t\right)  \\ 
\quad \qquad -f_{y}\left( t,y\left( t\right) ,Y\left( t\right) ,z\left(
t\right) ,Z\left( t\right) ,u\left( t\right) \right) q\left( t\right)  \\ 
\quad \qquad +G_{y}\left( t,y\left( t\right) ,Y\left( t\right) ,z\left(
t\right) ,Z\left( t\right) ,u\left( t\right) \right) k\left( t\right)  \\ 
\quad \qquad -g_{y}\left( t,y\left( t\right) ,Y\left( t\right) ,z\left(
t\right) ,Z\left( t\right) ,u\left( t\right) \right) h\left( t\right)  \\ 
\quad \qquad -l_{y}\left( t,y\left( t\right) ,Y\left( t\right) ,z\left(
t\right) ,Z\left( t\right) ,u\left( t\right) \right) ]\text{d}t \\ 
\quad \qquad +[F_{z}\left( t,y\left( t\right) ,Y\left( t\right) ,z\left(
t\right) ,Z\left( t\right) ,u\left( t\right) \right) p\left( t\right)  \\ 
\quad \qquad -f_{z}\left( t,y\left( t\right) ,Y\left( t\right) ,z\left(
t\right) ,Z\left( t\right) ,u\left( t\right) \right) q\left( t\right)  \\ 
\quad \qquad +G_{z}\left( t,y\left( t\right) ,Y\left( t\right) ,z\left(
t\right) ,Z\left( t\right) ,u\left( t\right) \right) k\left( t\right)  \\ 
\quad \qquad -g_{z}\left( t,y\left( t\right) ,Y\left( t\right) ,z\left(
t\right) ,Z\left( t\right) ,u\left( t\right) \right) h\left( t\right)  \\ 
\quad \qquad -l_{z}\left( t,y\left( t\right) ,Y\left( t\right) ,z\left(
t\right) ,Z\left( t\right) ,u\left( t\right) \right) ]\text{d}\overleftarrow{%
B_{t}}+h_{t}\text{d}\overrightarrow{W_{t}}, \\ 
p\left( 0\right) =-\gamma _{Y}\left( Y\left( 0\right) \right) ,\quad q\left(
T\right) =-\varphi _{y}\left( y\left( T\right) \right) P\left( T\right)
+\Phi _{y}\left( y\left( T\right) \right) .%
\end{array}%
\right.   \tag{3.2}
\end{equation}%
It is easy to check that FBDSDE (3.2) satisfies (H1), (H2) and (H$^{^{\prime
}}$3), so it has a unique solution 
\begin{equation*}
\left( p\left( t\right) ,q\left( t\right) ,k\left( t\right) ,h\left(
t\right) \right) \in M^{2}\left( 0,T;\mathbb{R}^{n}\mathbb{\times R}^{n}%
\mathbb{\times R}^{n\times l}\mathbb{\times R}^{n\times d}\right) .
\end{equation*}

We define the Hamiltonian function $H$ as follows: 
\begin{eqnarray*}
&&H\left( t,y\left( t\right) ,Y\left( t\right) ,z\left( t\right) ,Z\left(
t\right) ,v\left( t\right) ,p\left( t\right) ,q\left( t\right) ,k\left(
t\right) ,h\left( t\right) \right) \\
&\doteq &\left\langle q\left( t\right) ,f\left( t,y\left( t\right) ,Y\left(
t\right) ,z\left( t\right) ,Z\left( t\right) ,v\left( t\right) \right)
\right\rangle \\
&&-\left\langle p\left( t\right) ,F\left( t,y\left( t\right) ,Y\left(
t\right) ,z\left( t\right) ,Z\left( t\right) ,v\left( t\right) \right)
\right\rangle \\
&&-\left\langle k\left( t\right) ,G\left( t,y\left( t\right) ,Y\left(
t\right) ,z\left( t\right) ,Z\left( t\right) ,v\left( t\right) \right)
\right\rangle \\
&&+\left\langle h\left( t\right) ,g\left( t,y\left( t\right) ,Y\left(
t\right) ,z\left( t\right) ,Z\left( t\right) ,v\left( t\right) \right)
\right\rangle \\
&&+l\left( t,y\left( t\right) ,Y\left( t\right) ,z\left( t\right) ,Z\left(
t\right) ,v\left( t\right) \right) .
\end{eqnarray*}%
\begin{equation}
\tag{3.3}
\end{equation}%
FBDSDEs (3.2) can be rewritten as 
\begin{equation}
\left\{ 
\begin{array}{c}
\text{d}p\left( t\right) =-H_{Y}\text{d}t-H_{Z}\text{d}\overrightarrow{W_{t}}%
-k\left( t\right) \text{d}\overleftarrow{B_{t}}, \\ 
\text{d}q\left( t\right) =-H_{y}\text{d}t-H_{z}\text{d}\overleftarrow{B_{t}}%
+h\left( t\right) \text{d}\overrightarrow{W_{t}}, \\ 
q\left( T\right) =-\varphi _{y}\left( y\left( T\right) \right) p\left(
T\right) +\Phi _{y}\left( y\left( T\right) \right) , \\ 
p\left( 0\right) =-\gamma _{Y}\left( Y\left( 0\right) \right) ,\quad \quad
0\leq t\leq T,%
\end{array}%
\right.  \tag{3.4}
\end{equation}%
where $H_{\beta }=H_{\beta }\left( t,y\left( t\right) ,Y\left( t\right)
,z\left( t\right) ,Z\left( t\right) ,u\left( t\right) ,p\left( t\right)
,q\left( t\right) ,k\left( t\right) ,h\left( t\right) \right) ,$ $\beta =y,$ 
$Y,$ $z,$ $Z,$ respectively. At last, we can claim the first major result in
this paper.

\begin{theorem}[\textbf{Necessary maximum principle}]
Let $u\left( \cdot \right) $ be an optimal control and let $\left( y\left(
\cdot \right) ,Y\left( \cdot \right) ,z\left( \cdot \right) ,Z\left( \cdot
\right) \right) $ be the corresponding trajectory. Then we have 
\begin{eqnarray*}
\left\langle H_{v}\left( t,y\left( t\right) ,Y\left( t\right) ,z\left(
t\right) ,Z\left( t\right) ,u\left( t\right) ,p\left( t\right) ,q\left(
t\right) ,k\left( t\right) ,h\left( t\right) \right) ,v-u\left( t\right)
\right\rangle &\geq &0, \\
\text{ a.e., a.s. }t &\in &\left[ 0,T\right] ,\text{ }\forall v\in \mathcal{U%
},
\end{eqnarray*}%
\begin{equation}
\tag{3.5}
\end{equation}%
where $\left( p\left( t\right) ,q\left( t\right) ,k\left( t\right) ,h\left(
t\right) \right) $ is the solution of the adjoint equation (3.2).
\end{theorem}

\begin{proof}
Applying Itô formula to $\left\langle y^{1}\left( t\right) ,q\left( t\right)
\right\rangle +\left\langle Y^{1}\left( t\right) ,p\left( t\right)
\right\rangle $ on $\left[ 0,T\right] ,$ we have
\begin{eqnarray*}
&&\mathbb{E}\left[ \left\langle y^{1}\left( T\right) ,q\left( T\right)
\right\rangle +\left\langle Y^{1}\left( T\right) ,p\left( T\right)
\right\rangle -\left\langle y^{1}\left( 0\right) ,q\left( 0\right)
\right\rangle -\left\langle Y^{1}\left( 0\right) ,p\left( 0\right)
\right\rangle \right] \\
&&+\mathbb{E}\int_{0}^{T}[l_{y}\left( t,y\left( t\right) ,Y\left( t\right)
,z\left( t\right) ,Z\left( t\right) ,u\left( t\right) \right) y^{1}\left(
t\right) \\
&&+l_{Y}\left( t,y\left( t\right) ,Y\left( t\right) ,z\left( t\right)
,Z\left( t\right) ,u\left( t\right) \right) Y^{1}\left( t\right) \\
&&+l_{z}\left( t,y\left( t\right) ,Y\left( t\right) ,z\left( t\right)
,Z\left( t\right) ,u\left( t\right) \right) z^{1}\left( t\right) \\
&&+l_{Z}\left( t,y\left( t\right) ,Y\left( t\right) ,z\left( t\right)
,Z\left( t\right) ,u\left( t\right) \right) Z^{1}\left( t\right) \\
&&+l_{v}\left( t,y\left( t\right) ,Y\left( t\right) ,z\left( t\right)
,Z\left( t\right) ,u\left( t\right) \right) v\left( t\right) ]\text{d}t \\
&=&\mathbb{E}\int_{0}^{T}[\left\langle q\left( t\right) ,f_{v}\left(
t,y\left( t\right) ,Y\left( t\right) ,z\left( t\right) ,Z\left( t\right)
,u\left( t\right) \right) v\left( t\right) \right\rangle \\
&&-\left\langle p\left( t\right) ,F_{v}\left( t,y\left( t\right) ,Y\left(
t\right) ,z\left( t\right) ,Z\left( t\right) ,u\left( t\right) \right)
v\left( t\right) \right\rangle \\
&&-\left\langle k\left( t\right) ,G_{v}\left( t,y\left( t\right) ,Y\left(
t\right) ,z\left( t\right) ,Z\left( t\right) ,u\left( t\right) \right)
v\left( t\right) \right\rangle \\
&&+\left\langle h\left( t\right) ,g_{v}\left( t,y\left( t\right) ,Y\left(
t\right) ,z\left( t\right) ,Z\left( t\right) ,u\left( t\right) \right)
v\left( t\right) \right\rangle \\
&&+\left\langle v\left( t\right) ,l_{v}\left( t,y\left( t\right) ,Y\left(
t\right) ,z\left( t\right) ,Z\left( t\right) ,u\left( t\right) \right)
\right\rangle ]\text{d}t.
\end{eqnarray*}%
From the variational inequality in Lemma 5 and noting (3.3), for any $%
v\left( \cdot \right) \in \mathcal{U}_{ad}$ such that $u\left( \cdot \right)
+v\left( \cdot \right) \in \mathcal{U}_{ad},$ we have
\begin{equation*}
\mathbb{E}\int_{0}^{T}\left\langle H_{v}\left( t,y\left( t\right) ,Y\left(
t\right) ,z\left( t\right) ,Z\left( t\right) ,u\left( t\right) ,p\left(
t\right) ,q\left( t\right) ,k\left( t\right) ,h\left( t\right) \right)
,v\left( t\right) \right\rangle \text{d}t\geq 0.
\end{equation*}%
For $\forall v\in \mathcal{U},$ we set
\begin{equation*}
v\left( t\right) =\left\{
\begin{array}{l}
0,\qquad t\in \lbrack 0,t), \\
v,\qquad t\in \lbrack t,t+\varepsilon ), \\
0,\qquad t\in \left[ t+\varepsilon ,T\right] .%
\end{array}%
\right.
\end{equation*}%
Then we have
\begin{equation*}
\mathbb{E}\int_{t}^{t+\varepsilon }\left\langle H_{v}\left( t,y\left(
t\right) ,Y\left( t\right) ,z\left( t\right) ,Z\left( t\right) ,u\left(
t\right) ,p\left( t\right) ,q\left( t\right) ,k\left( t\right) ,h\left(
t\right) \right) ,v\right\rangle \text{d}t\geq 0.
\end{equation*}%
Notice the fact that
\begin{equation*}
\mathbb{E}\int_{t}^{t+\varepsilon }\left\langle H_{v}\left( t,y\left(
t\right) ,Y\left( t\right) ,z\left( t\right) ,Z\left( t\right) ,u\left(
t\right) ,p\left( t\right) ,q\left( t\right) ,k\left( t\right) ,h\left(
t\right) \right) ,u\left( t\right) \right\rangle \text{d}t=0.
\end{equation*}%
Differentiating with respect to $\varepsilon $ at $\varepsilon =0$ gives
\begin{eqnarray*}
\mathbb{E}\left\langle H_{v}\left( t,,y\left( t\right) ,Y\left( t\right)
,z\left( t\right) ,Z\left( t\right) ,u\left( t\right) ,p\left( t\right)
,q\left( t\right) ,k\left( t\right) ,h\left( t\right) \right) ,v-u\left(
t\right) \right\rangle &\geq &0, \\
\text{a.e., }t &\in &\left[ 0,T\right] .
\end{eqnarray*}%
The proof is completed.
\end{proof}

\section{A sufficient maximum principle for optimal controls of
forward-backward doubly stochastic control systems}

In this section, we investigate a sufficient maximum principle for the
optimal control problem stated in Section 2. For simplicity of notations, we
use the subscript label.

\begin{theorem}[\textbf{Sufficient maximum principle}]
Let $\left( \tilde{u}_{t};\tilde{y}_{t},\tilde{Y}_{t},\tilde{z}_{t},\tilde{Z}%
_{t}\right) _{0\leq t\leq T}$ be a quintuple and suppose there exist a
solution $\left( \tilde{p}_{t},\tilde{q}_{t},\tilde{k}_{t},\tilde{h}%
_{t}\right) _{0\leq t\leq T}$ of the corresponding adjoint forward-backward
doubly stochastic equation (3.2) such that for arbitrary admissible control $%
v\left( \cdot \right) \in \mathcal{U}_{ad},$ we have 
\begin{equation}
\mathbb{E}\int_{0}^{T}\left\langle \tilde{k}_{t},\left( Y_{t}-\tilde{Y}%
_{t}\right) \right\rangle ^{2}\text{d}t<\infty ,  \tag{4.1}
\end{equation}%
\begin{equation}
\mathbb{E}\int_{0}^{T}\left\langle \tilde{p}_{t},\left( Z_{t}-\tilde{Z}%
_{t}\right) \right\rangle ^{2}\text{d}t<\infty ,  \tag{4.2}
\end{equation}%
\begin{equation}
\mathbb{E}\int_{0}^{T}\left\langle \tilde{h}_{t},\left( y_{t}-\tilde{y}%
_{t}\right) \right\rangle ^{2}\text{d}t<\infty ,  \tag{4.3}
\end{equation}%
\begin{equation}
\mathbb{E}\int_{0}^{T}\left\langle \tilde{q}_{t},\left( z_{t}-\tilde{z}%
_{t}\right) \right\rangle ^{2}\text{d}t<\infty ,  \tag{4.4}
\end{equation}%
\begin{equation}
\mathbb{E}\int_{0}^{T}\left\langle \left( Y_{t}-\tilde{Y}_{t}\right)
,H_{Z}\left( t,\tilde{y}_{t},\tilde{Y}_{t},\tilde{z}_{t},\tilde{Z}_{t},%
\tilde{u}_{t},\tilde{p}_{t},\tilde{q}_{t},\tilde{k}_{t},\tilde{h}_{t}\right)
\right\rangle ^{2}\text{d}t<\infty ,  \tag{4.5}
\end{equation}%
\begin{equation}
\mathbb{E}\int_{0}^{T}\left\langle \tilde{p}_{t},\left( G\left(
t,y_{t},Y_{t},z_{t},Z_{t}\right) -G\left( t,\tilde{y}_{t},\tilde{Y}_{t},%
\tilde{z}_{t},\tilde{Z}_{t}\right) \right) \right\rangle ^{2}\text{d}%
t<\infty ,  \tag{4.6}
\end{equation}%
\begin{equation}
\mathbb{E}\int_{0}^{T}\left\langle \left( y_{t}-\tilde{y}_{t}\right)
,H_{z}\left( t,\tilde{y}_{t},\tilde{Y}_{t},\tilde{z}_{t},\tilde{Z}_{t},%
\tilde{u}_{t},\tilde{p}_{t},\tilde{q}_{t},\tilde{k}_{t},\tilde{h}_{t}\right)
\right\rangle ^{2}\text{d}t<\infty ,  \tag{4.7}
\end{equation}%
\begin{equation}
\mathbb{E}\int_{0}^{T}\left\langle \tilde{q}_{t},\left( g\left(
t,y_{t},Y_{t},z_{t},Z_{t}\right) -g\left( t,\tilde{y}_{t},\tilde{Y}_{t},%
\tilde{z}_{t},\tilde{Z}_{t}\right) \right) \right\rangle ^{2}\text{d}%
t<\infty .  \tag{4.8}
\end{equation}%
Further, suppose that for all 
\begin{equation*}
H\left( t,y,Y,z,Z,v,\tilde{p}_{t},\tilde{q}_{t},\tilde{k}_{t},\tilde{h}%
_{t}\right) ,\text{ }t\in \left[ 0,T\right] ,\text{ }
\end{equation*}
is convex in $\left( y,Y,z,Z,v\right) ,$ and $\gamma \left( Y\right) $ is
convex in $Y$ and $\Phi $ is convex in $y,$ moreover the following
conditions holds 
\begin{equation}
\mathbb{E}\left[ H\left( t,\tilde{y}_{t},\tilde{Y}_{t},\tilde{z}_{t},\tilde{Z%
}_{t},\tilde{u}_{t},\tilde{p}_{t},\tilde{q}_{t},\tilde{k}_{t},\tilde{h}%
_{t}\right) \right] =\underset{v\in U}{\inf }\mathbb{E}\left[ H\left( t,%
\tilde{y}_{t},\tilde{Y}_{t},\tilde{z}_{t},\tilde{Z}_{t},v,\tilde{p}_{t},%
\tilde{q}_{t},\tilde{k}_{t},\tilde{h}_{t}\right) \right] .  \tag{4.9}
\end{equation}%
Then $\tilde{u}_{t}$ is an optimal control.
\end{theorem}

\begin{proof}
Let $\left( y_{t},Y_{t},z_{t},Z_{t},v_{t}\right) =\left( y_{t}^{\left(
v\right) },Y_{t}^{\left( v\right) },z_{t}^{\left( v\right) },Z_{t}^{\left(
v\right) },v_{t}\right) $ be an arbitrary quintuple satisfying the control
system (2.1). According to the definition of the cost function (2.2), we
have
\begin{eqnarray*}
J\left( v\left( \cdot \right) \right) -J\left( \tilde{u}\left( \cdot \right)
\right) &=&\mathbb{E}\int_{0}^{T}\left[ l\left(
t,y_{t},Y_{t},z_{t},Z_{t},v_{t}\right) -l\left( t,\tilde{y}_{t},\tilde{Y}%
_{t},\tilde{z}_{t},\tilde{Z}_{t},\tilde{u}_{t}\right) \right] \text{d}t \\
&&+\mathbb{E}\left[ \Phi \left( y_{T}\right) -\Phi \left( \tilde{y}%
_{T}\right) \right] +E\left[ \gamma \left( Y_{0}\right) -\gamma \left(
\tilde{Y}_{0}\right) \right] \\
&=&\mathbf{I}_{1}\mathbf{+I}_{2}\mathbf{+I}_{3},
\end{eqnarray*}%
where
\begin{eqnarray*}
\mathbf{I}_{1} &=&\mathbb{E}\int_{0}^{T}\left[ l\left(
t,y_{t},Y_{t},z_{t},Z_{t},v_{t}\right) -l\left( t,\tilde{y}_{t},\tilde{Y}%
_{t},\tilde{z}_{t},\tilde{Z}_{t},\tilde{u}_{t}\right) \right] \text{d}t, \\
\mathbf{I}_{2} &=&\mathbb{E}\left[ \Phi \left( y_{T}\right) -\Phi \left(
\tilde{y}_{T}\right) \right] , \\
\mathbf{I}_{3} &=&\mathbb{E}\left[ \gamma \left( Y_{0}\right) -\gamma \left(
\tilde{Y}_{0}\right) \right] .
\end{eqnarray*}%
Now applying Itô formula to $\left\langle \tilde{p}_{t},Y_{t}-\tilde{Y}%
_{t}\right\rangle +\left\langle \tilde{q}_{t},y_{t}-\tilde{y}%
_{t}\right\rangle $ on $\left[ 0,T\right] ,$ we get
\begin{eqnarray*}
&&\left\langle \tilde{p}_{T},Y_{T}-\tilde{Y}_{T}\right\rangle +\left\langle
\tilde{q}_{T},y_{T}-\tilde{y}_{T}\right\rangle -\left\langle \tilde{p}%
_{0},Y_{0}-\tilde{Y}_{0}\right\rangle -\left\langle \tilde{q}_{0},y_{0}-%
\tilde{y}_{0}\right\rangle \\
\ &=&\left\langle \Phi _{y}\left( \tilde{y}_{T}\right) ,y_{T}-\tilde{y}%
_{T}\right\rangle +\left\langle \gamma _{Y}\left( \tilde{Y}_{0}\right)
,Y_{0}-\tilde{Y}_{0}\right\rangle \\
\ &=&\int_{0}^{T}\left\langle \left( Z_{t}-\tilde{Z}_{t}\right) ,\left(
-H_{Z}\left( t,\tilde{y}_{t},\tilde{Y}_{t},\tilde{z}_{t},\tilde{Z}_{t},%
\tilde{u}_{t},\tilde{p}_{t},\tilde{q}_{t},\tilde{k}_{t},\tilde{h}_{t}\right)
\right) \right\rangle \text{d}t \\
&&\ -\int_{0}^{T}\left\langle \tilde{k}_{t},\left( G\left(
t,y_{t},Y_{t},z_{t},Z_{t},v_{t}\right) -G\left( t,\tilde{y}_{t},\tilde{Y}%
_{t},\tilde{z}_{t},\tilde{Z}_{t},\tilde{u}_{t}\right) \right) \right\rangle
\text{d}t \\
&&\ +\int_{0}^{T}\left\langle \left( z_{t}-\tilde{z}_{t}\right) ,\left(
-H_{z}\left( t,\tilde{y}_{t},\tilde{Y}_{t},\tilde{z}_{t},\tilde{Z}_{t},%
\tilde{u}_{t},\tilde{p}_{t},\tilde{q}_{t},\tilde{k}_{t},\tilde{h}_{t}\right)
\right) \right\rangle \text{d}t \\
&&\ +\int_{0}^{T}\left\langle \tilde{h}_{t},\left( g\left(
t,y_{t},Y_{t},z_{t},Z_{t},v_{t}\right) -g\left( t,\tilde{y}_{t},\tilde{Y}%
_{t},\tilde{z}_{t},\tilde{Z}_{t},\tilde{u}_{t}\right) \right) \right\rangle
\text{d}t \\
&&\ +\int_{0}^{T}\left\langle \left( Y_{t}-\tilde{Y}_{t}\right) ,\left(
-H_{Y}\left( t,\tilde{y}_{t},\tilde{Y}_{t},\tilde{z}_{t},\tilde{Z}_{t},%
\tilde{u}_{t},\tilde{p}_{t},\tilde{q}_{t},\tilde{k}_{t},\tilde{h}_{t}\right)
\right) \right\rangle \text{d}t \\
&&\ +\int_{0}^{T}\left\langle \left( Y_{t}-\tilde{Y}_{t}\right) ,\left(
-H_{Z}\left( t,\tilde{y}_{t},\tilde{Y}_{t},\tilde{z}_{t},\tilde{Z}_{t},%
\tilde{u}_{t},\tilde{p}_{t},\tilde{q}_{t},\tilde{k}_{t},\tilde{h}_{t}\right)
\right) \text{d}\overrightarrow{W_{t}}\right\rangle \\
&&\ -\int_{0}^{T}\left\langle \tilde{k}_{t},\left( Y_{t}-\tilde{Y}%
_{t}\right) \text{d}\overleftarrow{B_{t}}\right\rangle \\
&&\ -\int_{0}^{T}\left\langle \tilde{p}_{t},\left( F\left(
t,y_{t},Y_{t},z_{t},Z_{t},v_{t}\right) -F\left( t,\tilde{y}_{t},\tilde{Y}%
_{t},\tilde{z}_{t},\tilde{Z}_{t},\tilde{u}_{t}\right) \right) \right\rangle
\text{d}t \\
&&\ -\int_{0}^{T}\left\langle \tilde{p}_{t},\left( G\left(
t,y_{t},Y_{t},z_{t},Z_{t},v_{t}\right) -G\left( t,\tilde{y}_{t},\tilde{Y}%
_{t},\tilde{z}_{t},\tilde{Z}_{t},\tilde{u}_{t}\right) \right) \text{d}%
\overleftarrow{B_{t}}\right\rangle \\
&&\ +\int_{0}^{T}\left\langle \tilde{p}_{t},\left( Z_{t}-\tilde{Z}%
_{t}\right) \text{d}\overrightarrow{W_{t}}\right\rangle \\
&&\ +\int_{0}^{T}\left\langle \left( y_{t}-\tilde{y}_{t}\right) ,\left(
-H_{y}\left( t,\tilde{y}_{t},\tilde{Y}_{t},\tilde{z}_{t},\tilde{Z}_{t},%
\tilde{u}_{t},\tilde{p}_{t},\tilde{q}_{t},\tilde{k}_{t},\tilde{h}_{t}\right)
\right) \right\rangle \text{d}t \\
&&\ +\int_{0}^{T}\left\langle \left( y_{t}-\tilde{y}_{t}\right) ,\left(
-H_{z}\left( t,\tilde{y}_{t},\tilde{Y}_{t},\tilde{z}_{t},\tilde{Z}_{t},%
\tilde{u}_{t},\tilde{p}_{t},\tilde{q}_{t},\tilde{k}_{t},\tilde{h}_{t}\right)
\right) \text{d}\overleftarrow{B_{t}}\right\rangle \\
&&\ +\int_{0}^{T}\left\langle \left( y_{t}-\tilde{y}_{t}\right) ,\tilde{h}%
_{t}\text{d}W_{t}\right\rangle \\
&&\ +\int_{0}^{T}\tilde{q}_{t}\left( f\left(
t,y_{t},Y_{t},z_{t},Z_{t},v_{t}\right) -f\left( t,\tilde{y}_{t},\tilde{Y}%
_{t},\tilde{z}_{t},\tilde{Z}_{t},\tilde{u}_{t}\right) \right) \text{d}t \\
&&\ +\int_{0}^{T}\tilde{q}_{t}\left( g\left(
t,y_{t},Y_{t},z_{t},Z_{t},v_{t}\right) -g\left( t,\tilde{y}_{t},\tilde{Y}%
_{t},\tilde{z}_{t},\tilde{Z}_{t},\tilde{u}_{t}\right) \right) \text{d}W_{t}
\\
&&\ -\int_{0}^{T}\tilde{q}_{t}\left( z_{t}-\tilde{z}_{t}\right) \text{d}%
B_{t},
\end{eqnarray*}%
where we claim that
\begin{equation*}
\left\{
\begin{array}{l}
Y_{T}-\tilde{Y}_{T}=\varphi \left( y_{T}\right) -\varphi \left( \tilde{y}%
_{T}\right) =\varphi _{y}\left( \tilde{y}\left( T\right) \right) \left(
y\left( T\right) -\tilde{y}\left( T\right) \right) , \\
y_{0}-\tilde{y}_{0}=x_{0}-x_{0}=0, \\
\tilde{p}_{0}=-\gamma _{Y}\left( Y_{0}\right) , \\
\tilde{q}_{T}=\Phi _{y}\left( \tilde{y}_{T}\right) -\varphi _{y}\left(
\tilde{y}\left( T\right) \right) \tilde{p}\left( T\right) .%
\end{array}%
\right.
\end{equation*}%
By Davis inequality, under the conditions (4.1)-(4.8), we can ensure that
the stochastic integrals with respect to the Brownian motion have zero
expectations. Moreover, by virtue of convexity of $\Phi $ and $\gamma $, it
follows instantly that
\begin{eqnarray*}
\mathbf{I}_{2}+\mathbf{I}_{3} &=&\mathbb{E}\left[ \Phi \left( y_{T}\right)
-\Phi \left( \tilde{y}_{T}\right) \right] +\mathbb{E}\left[ \gamma \left(
Y_{0}\right) -\gamma \left( \tilde{Y}_{0}\right) \right] \\
&\geq &\mathbb{E}\left\langle \Phi _{y}\left( \tilde{y}_{T}\right) ,y_{T}-%
\tilde{y}_{T}\right\rangle +\mathbb{E}\left\langle \gamma _{Y}\left( \tilde{Y%
}_{0}\right) ,Y_{0}-\tilde{Y}_{0}\right\rangle \\
&=&-\mathbb{E}\int_{0}^{T}\left\langle \left( Y_{t}-\tilde{Y}_{t}\right)
,H_{Y}\left( t,\tilde{y}_{t},\tilde{Y}_{t},\tilde{z}_{t},\tilde{Z}_{t},%
\tilde{u}_{t},\tilde{p}_{t},\tilde{q}_{t},\tilde{k}_{t},\tilde{h}_{t}\right)
\right\rangle \text{d}t \\
&&-\mathbb{E}\int_{0}^{T}\left\langle \tilde{p}_{t},\left( F\left(
t,y_{t},Y_{t},z_{t},Z_{t},v_{t}\right) -F\left( t,\tilde{y}_{t},\tilde{Y}%
_{t},\tilde{z}_{t},\tilde{Z}_{t},\tilde{u}_{t}\right) \right) \right\rangle
\text{d}t \\
&&-\mathbb{E}\int_{0}^{T}\left\langle \left( y_{t}-\tilde{y}_{t}\right)
,H_{y}\left( t,\tilde{y}_{t},\tilde{Y}_{t},\tilde{z}_{t},\tilde{Z}_{t},%
\tilde{u}_{t},\tilde{p}_{t},\tilde{q}_{t},\tilde{k}_{t},\tilde{h}_{t}\right)
\right\rangle \text{d}t \\
&&+\mathbb{E}\int_{0}^{T}\left\langle \tilde{q}_{t},\left( g\left(
t,y_{t},Y_{t},z_{t},Z_{t},v_{t}\right) -g\left( t,\tilde{y}_{t},\tilde{Y}%
_{t},\tilde{z}_{t},\tilde{Z}_{t},\tilde{u}_{t}\right) \right) \right\rangle
\text{d}t \\
&&-\mathbb{E}\int_{0}^{T}\left\langle \left( Z_{t}-\tilde{Z}_{t}\right)
,H_{Z}\left( t,\tilde{y}_{t},\tilde{Y}_{t},\tilde{z}_{t},\tilde{Z}_{t},%
\tilde{u}_{t},\tilde{p}_{t},\tilde{q}_{t},\tilde{k}_{t},\tilde{h}_{t}\right)
\right\rangle \text{d}t \\
&&-\mathbb{E}\int_{0}^{T}\left\langle \tilde{k}_{t},\left( G\left(
t,y_{t},Y_{t},z_{t},Z_{t},v_{t}\right) -G\left( t,\tilde{y}_{t},\tilde{Y}%
_{t},\tilde{z}_{t},\tilde{Z}_{t},\tilde{u}_{t}\right) \right) \right\rangle
\text{d}t \\
&&-\mathbb{E}\int_{0}^{T}\left\langle \left( z_{t}-\tilde{z}_{t}\right)
,H_{z}\left( t,\tilde{y}_{t},\tilde{Y}_{t},\tilde{z}_{t},\tilde{Z}_{t},%
\tilde{u}_{t},\tilde{p}_{t},\tilde{q}_{t},\tilde{k}_{t},\tilde{h}_{t}\right)
\right\rangle \text{d}t \\
&&+\mathbb{E}\int_{0}^{T}\left\langle \tilde{h}_{t},\left( g\left(
t,y_{t},Y_{t},z_{t},Z_{t},v_{t}\right) -g\left( t,\tilde{y}_{t},\tilde{Y}%
_{t},\tilde{z}_{t},\tilde{Z}_{t},\tilde{u}_{t}\right) \right) \right\rangle
\text{d}t \\
&=&-\Xi _{1}+\Xi _{2}+\Xi _{3}+\Xi _{4}+\Xi _{5},
\end{eqnarray*}%
where
\begin{eqnarray*}
\Xi _{1} &=&\mathbb{E}\int_{0}^{T}\left\langle H_{y}\left( t,\tilde{y}_{t},%
\tilde{Y}_{t},\tilde{z}_{t},\tilde{Z}_{t},\tilde{u}_{t},\tilde{p}_{t},\tilde{%
q}_{t},\tilde{k}_{t},\tilde{h}_{t}\right) ,\left( y_{t}-\tilde{y}_{t}\right)
\right\rangle \text{d}t \\
&&+\mathbb{E}\int_{0}^{T}\left\langle H_{Y}\left( t,\tilde{y}_{t},\tilde{Y}%
_{t},\tilde{z}_{t},\tilde{Z}_{t},\tilde{u}_{t},\tilde{p}_{t},\tilde{q}_{t},%
\tilde{k}_{t},\tilde{h}_{t}\right) ,\left( Y_{t}-\tilde{Y}_{t}\right)
\right\rangle \text{d}t \\
&&+\mathbb{E}\int_{0}^{T}\left\langle H_{z},\left( t,\tilde{y}_{t},\tilde{Y}%
_{t},\tilde{z}_{t},\tilde{Z}_{t},\tilde{u}_{t},\tilde{p}_{t},\tilde{q}_{t},%
\tilde{k}_{t},\tilde{h}_{t}\right) ,\left( z_{t}-\tilde{z}_{t}\right)
\right\rangle \text{d}t \\
&&+\mathbb{E}\int_{0}^{T}\left\langle H_{Z}\left( t,\tilde{y}_{t},\tilde{Y}%
_{t},\tilde{z}_{t},\tilde{Z}_{t},\tilde{u}_{t},\tilde{p}_{t},\tilde{q}_{t},%
\tilde{k}_{t},\tilde{h}_{t}\right) ,\left( Z_{t}-\tilde{Z}_{t}\right)
\right\rangle \text{d}t \\
\Xi _{2} &=&-\mathbb{E}\int_{0}^{T}\left\langle \tilde{p}_{t},F\left(
t,y_{t},Y_{t},z_{t},Z_{t},v_{t}\right) -F\left( t,\tilde{y}_{t},\tilde{Y}%
_{t},\tilde{z}_{t},\tilde{Z}_{t},\tilde{u}_{t}\right) \right\rangle \text{d}t
\\
\Xi _{3} &=&\mathbb{E}\int_{0}^{T}\left\langle \tilde{q}_{t},g\left(
t,y_{t},Y_{t},z_{t},Z_{t},v_{t}\right) -g\left( t,\tilde{y}_{t},\tilde{Y}%
_{t},\tilde{z}_{t},\tilde{Z}_{t},\tilde{u}_{t}\right) \right\rangle \text{d}t
\\
\Xi _{4} &=&-\mathbb{E}\int_{0}^{T}\left\langle \tilde{k}_{t},G\left(
t,y_{t},Y_{t},z_{t},Z_{t},v_{t}\right) -G\left( t,\tilde{y}_{t},\tilde{Y}%
_{t},\tilde{z}_{t},\tilde{Z}_{t},\tilde{u}_{t}\right) \right\rangle \text{d}t
\\
\Xi _{5} &=&\mathbb{E}\int_{0}^{T}\left\langle \tilde{h}_{t},g\left(
t,y_{t},Y_{t},z_{t},Z_{t},v_{t}\right) -g\left( t,\tilde{y}_{t},\tilde{Y}%
_{t},\tilde{z}_{t},\tilde{Z}_{t},\tilde{u}_{t}\right) \right\rangle \text{d}%
t.
\end{eqnarray*}%
Noting the definition of $H$ and $\mathbf{I}_{1},$ we have
\begin{eqnarray*}
\mathbf{I}_{1} &=&\mathbb{E}\int_{0}^{T}\left[ l\left(
t,y_{t},Y_{t},z_{t},Z_{t},v_{t}\right) -l\left( t,\tilde{y}_{t},\tilde{Y}%
_{t},\tilde{z}_{t},\tilde{Z}_{t},\tilde{u}_{t}\right) \right] \text{d}t \\
&=&\mathbb{E}\int_{0}^{T}\left[ H\left( t,y_{t},Y_{t},z_{t},Z_{t},v_{t},%
\tilde{p}_{t},\tilde{q}_{t},\tilde{k}_{t},\tilde{h}_{t}\right) -H\left( t,%
\tilde{y}_{t},\tilde{Y}_{t},\tilde{z}_{t},\tilde{Z}_{t},\tilde{u}_{t},\tilde{%
p}_{t},\tilde{q}_{t},\tilde{k}_{t},\tilde{h}_{t}\right) \right] \text{d}t \\
&&-\mathbb{E}\int_{0}^{T}\left[ \left\langle \tilde{q}_{t},f\left(
t,y_{t},Y_{t},z_{t},Z_{t},v_{t}\right) -f\left( t,\tilde{y}_{t},\tilde{Y}%
_{t},\tilde{z}_{t},\tilde{Z}_{t},\tilde{u}_{t}\right) \right\rangle \right]
\text{d}t \\
&&+\mathbb{E}\int_{0}^{T}\left[ \left\langle \tilde{p}_{t},F\left(
t,y_{t},Y_{t},z_{t},Z_{t},v_{t}\right) -F\left( t,\tilde{y}_{t},\tilde{Y}%
_{t},\tilde{z}_{t},\tilde{Z}_{t},\tilde{u}_{t}\right) \right\rangle \right]
\text{d}t \\
&&+\mathbb{E}\int_{0}^{T}\left[ \left\langle \tilde{k}_{t},G\left(
t,y_{t},Y_{t},z_{t},Z_{t},v_{t}\right) -G\left( t,\tilde{y}_{t},\tilde{Y}%
_{t},\tilde{z}_{t},\tilde{Z}_{t},\tilde{u}_{t}\right) \right\rangle \right]
\text{d}t \\
&&-\mathbb{E}\int_{0}^{T}\left[ \left\langle \tilde{h}_{t},g\left(
t,y_{t},Y_{t},z_{t},Z_{t},v_{t}\right) -g\left( t,\tilde{y}_{t},\tilde{Y}%
_{t},\tilde{z}_{t},\tilde{Z}_{t},\tilde{u}_{t}\right) \right\rangle \right]
\text{d}t \\
&=&\Xi _{6}-\Xi _{2}-\Xi _{3}-\Xi _{4}-\Xi _{5},
\end{eqnarray*}%
where
\begin{equation*}
\Xi _{6}=\mathbb{E}\int_{0}^{T}\left[ H\left(
t,y_{t},Y_{t},z_{t},Z_{t},v_{t},\tilde{p}_{t},\tilde{q}_{t},\tilde{k}_{t},%
\tilde{h}_{t}\right) -H\left( t,\tilde{y}_{t},\tilde{Y}_{t},\tilde{z}_{t},%
\tilde{Z}_{t},\tilde{u}_{t},\tilde{p}_{t},\tilde{q}_{t},\tilde{k}_{t},\tilde{%
h}_{t}\right) \right] \text{d}t.
\end{equation*}%
On the one hand, by the virtue of convexity of $H\left( t,y,Y,z,Z,v,\tilde{p}%
_{t},\tilde{q}_{t},\tilde{k}_{t},\tilde{h}_{t}\right) $ with respect to $%
\left( y,Y,z,Z,v\right) ,$ we obtain
\begin{eqnarray*}
&&H\left( t,y_{t},Y_{t},z_{t},Z_{t},v_{t},\tilde{p}_{t},\tilde{q}_{t},\tilde{%
k}_{t},\tilde{h}_{t}\right) -H\left( t,\tilde{y}_{t},\tilde{Y}_{t},\tilde{z}%
_{t},\tilde{Z}_{t},\tilde{u}_{t},\tilde{p}_{t},\tilde{q}_{t},\tilde{k}_{t},%
\tilde{h}_{t}\right) \\
&\geq &H_{y}\left( t,\tilde{y}_{t},\tilde{Y}_{t},\tilde{z}_{t},\tilde{Z}_{t},%
\tilde{u}_{t},\tilde{p}_{t},\tilde{q}_{t},\tilde{k}_{t},\tilde{h}_{t}\right)
\left( y_{t}-\tilde{y}_{t}\right) \\
&&+H_{Y}\left( t,\tilde{y}_{t},\tilde{Y}_{t},\tilde{z}_{t},\tilde{Z}_{t},%
\tilde{u}_{t},\tilde{p}_{t},\tilde{q}_{t},\tilde{k}_{t},\tilde{h}_{t}\right)
\left( Y_{t}-\tilde{Y}_{t}\right) \\
&&+H_{z}\left( t,\tilde{y}_{t},\tilde{Y}_{t},\tilde{z}_{t},\tilde{Z}_{t},%
\tilde{u}_{t},\tilde{p}_{t},\tilde{q}_{t},\tilde{k}_{t},\tilde{h}_{t}\right)
\left( z_{t}-\tilde{z}_{t}\right) \\
&&+H_{Z}\left( t,\tilde{y}_{t},\tilde{Y}_{t},\tilde{z}_{t},\tilde{Z}_{t},%
\tilde{u}_{t},\tilde{p}_{t},\tilde{q}_{t},\tilde{k}_{t},\tilde{h}_{t}\right)
\left( Z_{t}-\tilde{Z}_{t}\right) \\
&&+H_{u}\left( t,\tilde{y}_{t},\tilde{Y}_{t},\tilde{z}_{t},\tilde{Z}_{t},%
\tilde{u}_{t},\tilde{p}_{t},\tilde{q}_{t},\tilde{k}_{t},\tilde{h}_{t}\right)
\left( v_{t}-\tilde{u}_{t}\right)
\end{eqnarray*}%
\begin{equation}
\tag{4.10}
\end{equation}%
On the other hand, we know
\begin{equation*}
\mathbb{E}\left[ H_{u}\left( t,\tilde{y}_{t},\tilde{Y}_{t},\tilde{z}_{t},%
\tilde{Z}_{t},\tilde{u}_{t},\tilde{p}_{t},\tilde{q}_{t},\tilde{k}_{t},\tilde{%
h}_{t}\right) \left( v_{t}-\tilde{u}_{t}\right) \right] \geq 0.
\end{equation*}%
Consequently, associating with (4.10), we claim that
\begin{eqnarray*}
\Xi _{6} &=&\mathbb{E}\int_{0}^{T}\left[ H\left(
t,y_{t},Y_{t},z_{t},Z_{t},v_{t},\tilde{p}_{t},\tilde{q}_{t},\tilde{k}_{t},%
\tilde{h}_{t}\right) -H\left( t,\tilde{y}_{t},\tilde{Y}_{t},\tilde{z}_{t},%
\tilde{Z}_{t},\tilde{u}_{t},\tilde{p}_{t},\tilde{q}_{t},\tilde{k}_{t},\tilde{%
h}_{t}\right) \right] \text{d}t \\
&\geq &\mathbb{E}\int_{0}^{T}H_{y}\left( t,\tilde{y}_{t},\tilde{Y}_{t},%
\tilde{z}_{t},\tilde{Z}_{t},\tilde{u}_{t},\tilde{p}_{t},\tilde{q}_{t},\tilde{%
k}_{t},\tilde{h}_{t}\right) \left( y_{t}-\tilde{y}_{t}\right) \text{d}t \\
&&+\mathbb{E}\int_{0}^{T}H_{Y}\left( t,\tilde{y}_{t},\tilde{Y}_{t},\tilde{z}%
_{t},\tilde{Z}_{t},\tilde{u}_{t},\tilde{p}_{t},\tilde{q}_{t},\tilde{k}_{t},%
\tilde{h}_{t}\right) \left( Y_{t}-\tilde{Y}_{t}\right) \text{d}t \\
&&+\mathbb{E}\int_{0}^{T}H_{z}\left( t,\tilde{y}_{t},\tilde{Y}_{t},\tilde{z}%
_{t},\tilde{Z}_{t},\tilde{u}_{t},\tilde{p}_{t},\tilde{q}_{t},\tilde{k}_{t},%
\tilde{h}_{t}\right) \left( z_{t}-\tilde{z}_{t}\right) \text{d}t \\
&&+\mathbb{E}\int_{0}^{T}H_{Z}\left( t,\tilde{y}_{t},\tilde{Y}_{t},\tilde{z}%
_{t},\tilde{Z}_{t},\tilde{u}_{t},\tilde{p}_{t},\tilde{q}_{t},\tilde{k}_{t},%
\tilde{h}_{t}\right) \left( Z_{t}-\tilde{Z}_{t}\right) \text{d}t \\
&=&\Xi _{1}.
\end{eqnarray*}%
Then, it follows that
\begin{eqnarray*}
J\left( v\left( \cdot \right) \right) -J\left( u\left( \cdot \right) \right)
&=&\mathbf{I}_{1}\mathbf{+I}_{2}\mathbf{+I}_{3} \\
&=&\Xi _{6}-\Xi _{2}-\Xi _{3}-\Xi _{4}-\Xi _{5} \\
&&-\Xi _{1}+\Xi _{2}+\Xi _{3}+\Xi _{4}+\Xi _{5} \\
&\geq &\Xi _{1}-\Xi _{2}-\Xi _{3}-\Xi _{4}-\Xi _{5} \\
&&-\Xi _{1}+\Xi _{2}+\Xi _{3}+\Xi _{4}+\Xi _{5} \\
&=&0.
\end{eqnarray*}%
Since $v\left( \cdot \right) \in \mathcal{U}_{ad}$ is arbitrary, we say that
$\tilde{u}\left( \cdot \right) $ is an optimal control. The proof is
completed.
\end{proof}

\section{Applications to optimal control problems of stochastic partial
differential equations}

In this section, we will give necessary and sufficient maximum principles
for optimal control of SPDEs. Let us first give some notations from [15].
For convenience, all the variables in this section are one-dimensional. It
is necessary to point out that all the results in this section can be
extended to multi-dimensional cases, but we use the notations in general
case. From now on $C^{k}\left( \mathbb{R}^{n};\mathbb{R}^{n}\right) ,$ $%
C_{l,b}^{k}\left( \mathbb{R}^{n};\mathbb{R}^{n}\right) ,$ $C_{p}^{k}\left( 
\mathbb{R}^{n};\mathbb{R}^{n}\right) $ will denote respectively the set of
functions of class $C^{k}$ from $\mathbb{R}^{n}$ into $\mathbb{R}^{n}$, the
set of those functions of class $C^{k}$ whose partial derivatives of order
less than or equal to $k$ are bounded (and hence the function itself grows
at most linearly at infinity), and the set of those functions of class $C^{k}
$ which, together with all their partial derivatives of order less than or
equal to $k$, grow at most like a polynomial function of the variable $x$ at
infinity. We consider the following quasilinear SPDEs with control variable: 
\begin{equation}
\left\{ 
\begin{array}{c}
u\left( t,x\right) =\varphi \left( x\right) +\int_{t}^{T}\left[ \mathcal{L}%
^{v}u\left( s,x\right) +F\left( s,x,u\left( s,x\right) ,\left( \nabla
u\sigma \right) \left( s,x,u\right) ,v\left( s\right) \right) \right] \text{d%
}s \\ 
+\int_{t}^{T}G\left( x,u\left( s,x\right) ,\left( \nabla u\sigma \right)
\left( s,x,u\right) ,v\left( s\right) \right) \text{d}\overleftarrow{B_{s}}%
,\quad 0\leq t\leq T,%
\end{array}%
\right.   \tag{5.1}
\end{equation}%
where $u:\left[ 0,T\right] \times \mathbb{R}^{d}\mathbb{\rightarrow R}^{k}$
and $\nabla u\left( s,x\right) $ denotes the first order derivative of $%
u\left( s,x\right) $ with respect to $x$, and 
\begin{equation*}
\mathcal{L}^{v}u=\left( 
\begin{array}{c}
L^{v}u_{1} \\ 
\vdots  \\ 
L^{v}u_{k}%
\end{array}%
\right) ,
\end{equation*}%
with 
\begin{equation*}
L\phi \left( x\right) =\frac{1}{2}\sum_{i,j=1}^{d}\left( gg^{\ast }\right)
_{ij}\left( x,v\right) \frac{\partial ^{2}\phi \left( x\right) }{\partial
x_{i}\partial x_{j}}+\sum_{i=1}^{d}f_{i}\left( x,v\right) \frac{\partial
\phi \left( x\right) }{\partial x_{i}}.
\end{equation*}%
and 
\begin{eqnarray*}
f &:&\mathbb{R}^{n}\mathbb{\times R}^{k}\mathbb{\rightarrow R}^{n}, \\
g &:&\mathbb{R}^{n}\mathbb{\times R}^{k}\mathbb{\rightarrow R}^{n}\mathbf{,}
\\
F &:&\left[ 0,T\right] \times \mathbb{R}^{n}\mathbb{\times R\times R}^{d}%
\mathbb{\times R}^{k}\mathbb{\rightarrow R}\mathbf{,} \\
G &:&\left[ 0,T\right] \times \mathbb{R}^{n}\mathbb{\times R\times R}^{d}%
\mathbb{\times R}^{k}\mathbb{\rightarrow R}\mathbf{,} \\
\varphi  &:&\mathbb{R}^{n}\rightarrow \mathbb{R}.
\end{eqnarray*}

In order to assure the existence and uniqueness of solutions for (5.1) and
(5.3) below, we give the following assumptions for sake of completeness (see
[15] for more details).

\begin{enumerate}
\item[\textbf{(A1)}] Assume that%
\begin{equation*}
\left\{ 
\begin{array}{l}
f\in C_{l,b}^{3}\left( \mathbb{R}^{n}\mathbb{\times R}^{k};\mathbb{R}%
^{n}\right) ,\quad g\in C_{l,b}^{3}\left( \mathbb{R}^{n}\mathbb{\times R}%
^{k};\mathbb{R}^{n\times d}\right) ,\text{\quad }\varphi \in C_{p}^{3}\left( 
\mathbb{R}^{n};\mathbb{R}\right) , \\ 
F\left( t,\cdot ,\cdot ,\cdot ,v\right) \in C_{l,b}^{3}\left( \mathbb{R}^{n}%
\mathbb{\times R\times R}^{d}\mathbb{\times R}^{k};\mathbb{R}\right) ,\quad
F\left( \cdot ,x,y,z,v\right) \in M^{2}\left( 0,T;\mathbb{R}\right) , \\ 
G\left( t,\cdot ,\cdot ,\cdot ,v\right) \in C_{l,b}^{3}\left( \mathbb{R}^{n}%
\mathbb{\times R\times R}^{d}\mathbb{\times R}^{k};\mathbb{R}\right) ,\quad
G\left( \cdot ,x,y,z,v\right) \in M^{2}\left( 0,T;\mathbb{R}\right) , \\ 
\forall t\in \left[ 0,T\right] \text{, }x\in \mathbb{R}^{n}\text{, }y\in 
\mathbb{R}\text{, }z\in \mathbb{R}^{d}\text{, }v\in \mathbb{R}^{k}.%
\end{array}%
\right. 
\end{equation*}

\item[\textbf{(A2)}] Assume that (H1), (H2) and (H3) hold.
\end{enumerate}

Let $\mathcal{U}_{ad}\footnote{%
Here 
\begin{equation*}
\mathcal{U}_{ad}\doteq \left\{ v\left( \cdot \right) \in M_{\mathcal{F}%
^{W}}^{2}\left( 0,T;\mathbb{R}^{k}\right) ;\text{ }v\left( t\right) \in 
\mathcal{U},\text{ }0\leq t\leq T,\text{ a.e., a.s.}\right\} .
\end{equation*}%
where$M_{\mathcal{F}^{W}}^{2}\left( 0,T;\mathbb{R}^{k}\right) $ denotes the
space of (class of $dP\otimes dt$ a.e equal) all $\left\{ \mathcal{F}%
_{t}^{W}\right\} $-measurable $k$-dimensional processes $\upsilon $ with
norm of $\parallel \upsilon \parallel _{M_{\mathcal{F}^{W}}^{2}}\doteq \left[
\mathbb{E}\int_{0}^{T}|\upsilon (s)|^{2}\text{d}s\right] ^{\frac{1}{2}%
}<\infty .$}$ be an admissible control set. The optimal control problem of
SPDEs (5.1) is to find an optimal control $v^{\ast }\left( \cdot \right) \in 
\mathcal{U}_{ad}$, such that 
\begin{equation*}
J\left( v^{\ast }\left( \cdot \right) \right) \doteq \underset{v\left( \cdot
\right) \in \mathcal{U}_{ad}}{\inf }J\left( v\left( \cdot \right) \right) ,
\end{equation*}%
where $J\left( \cdot \right) $ is the cost function as follows: 
\begin{equation}
J\left( v\left( \cdot \right) \right) =\mathbb{E}\left[ \int_{0}^{T}l\left(
s,\bar{x},u\left( s,\bar{x}\right) ,\left( \nabla ug\right) \left( s,\bar{x}%
,u\left( s,\bar{x}\right) \right) ,v\left( s\right) \right) _{\bar{x}%
=X^{0,x}\left( s\right) }\text{d}s+\gamma \left( u\left( 0,x\right) \right) %
\right] .  \tag{5.2}
\end{equation}%
Here we assume $l$ and $\gamma $ satisfy (H4) and $X^{0,x}\left( s\right) $
defined below. We can transform the optimal control problem of SPDEs (5.1)
into one of the following FBDSDEs with control variable $v\left( \cdot
\right) $: 
\begin{equation}
\left\{ 
\begin{array}{l}
X^{t,x}\left( s\right) =x+\int_{t}^{s}f\left( X^{t,x}\left( r\right)
,v\left( r\right) \right) \text{d}r+\int_{t}^{s}g\left( X^{t,x}\left(
r\right) ,v\left( r\right) \right) \text{d}\overrightarrow{W_{r}}, \\ 
Y^{t,x}\left( s\right) =\varphi \left( X^{t,x}\left( T\right) \right)
+\int_{s}^{T}F\left( r,X^{t,x}\left( r\right) ,Y^{t,x}\left( r\right)
,Z^{t,x}\left( r\right) ,v\left( r\right) \right) \text{d}r \\ 
\qquad \qquad +\int_{s}^{T}G\left( r,X^{t,x}\left( r\right) ,Y^{t,x}\left(
r\right) ,Z^{t,x}\left( r\right) ,v\left( r\right) \right) \text{d}%
\overleftarrow{B_{r}} \\ 
\qquad \qquad -\int_{s}^{T}Z^{t,x}\left( r\right) \text{d}\overrightarrow{%
W_{r}},\quad \quad 0\leq t\leq s\leq T,%
\end{array}%
\right.   \tag{5.3}
\end{equation}%
where $\left( X^{t,x}\left( \cdot \right) ,Y^{t,x}\left( \cdot \right)
,Z^{t,x}\left( \cdot \right) ,v\left( \cdot \right) \right) \in \mathbb{R}%
^{n}\mathbb{\times R\times R}^{d}\mathbb{\times R}^{k}$, $x\in \mathbb{R}^{n}
$. The corresponding optimal control problem of FBDSDEs (5.3) is to find an
optimal control $v^{\ast }\left( \cdot \right) \in \mathcal{U}_{ad}$, such
that 
\begin{equation*}
J\left( v^{\ast }\left( \cdot \right) \right) \doteq \underset{v\left( \cdot
\right) \in \mathcal{U}_{ad}}{\inf }J\left( v\left( \cdot \right) \right) ,
\end{equation*}%
where $J\left( v\left( \cdot \right) \right) $ is the cost function the same
as (5.2): 
\begin{equation*}
J\left( v\left( \cdot \right) \right) =\mathbb{E}\left[ \int_{0}^{T}l\left(
s,X\left( s\right) ,Y\left( s\right) ,Z\left( s\right) ,v\left( s\right)
\right) \text{d}s+\gamma \left( Y\left( 0\right) \right) \right] .
\end{equation*}%
Now we consider the following adjoint FBDSDEs involving the four unknown
processes $\left( p\left( t\right) ,q\left( t\right) ,k\left( t\right)
,h\left( t\right) \right) $: 
\begin{equation}
\left\{ 
\begin{array}{l}
\text{d}p\left( t\right) =\left( F_{Y}p\left( t\right) +G_{Y}k\left(
t\right) -l_{Y}\right) \text{d}t \\ 
\qquad \quad +\left( F_{Z}p\left( t\right) -G_{Z}k\left( t\right)
-l_{Z}\right) \text{d}\overrightarrow{W_{t}}-k\left( t\right) \text{d}%
\overleftarrow{B_{t}}, \\ 
\text{d}q\left( t\right) =\left( F_{X}p\left( t\right) -f_{X}q\left(
t\right) +G_{X}k\left( t\right) -g_{X}h\left( t\right) -l_{X}\right) \text{d}%
t+h\left( t\right) \text{d}\overrightarrow{W_{t}}, \\ 
p\left( 0\right) =-\gamma _{Y}\left( Y\left( 0\right) \right) ,\quad q\left(
T\right) =-\varphi _{X}\left( X\left( T\right) \right) p\left( T\right) ,%
\text{\qquad }0\leq t\leq T.%
\end{array}%
\right.   \tag{5.4}
\end{equation}%
It is easy to see that (5.4) satisfies (H1), (H2) and (H'3), so it is
uniquely solvable by virtue of Proposition 2. Therefore we know that (5.4)
has a unique solution $\left( p\left( \cdot \right) ,q\left( \cdot \right)
,k\left( \cdot \right) ,h\left( \cdot \right) \right) \in M^{2}\left( 0,T;%
\mathbb{R\times R}^{n}\mathbb{\times R}^{l}\mathbb{\times R}^{d}\right) $.
Define the Hamilton function as follows: 
\begin{eqnarray*}
\bar{H}\left( t,X,Y,Z,v,p,q,k,h\right)  &=&H\left(
t,X,Y,0,Z,v,p,q,k,h\right)  \\
&=&l\left( t,X,Y,Z,v\right) -k\cdot G\left( t,X,Y,Z,v\right)  \\
&&+q\cdot f\left( X,Y,v\right) -p\cdot F\left( t,X,Y,Z,v\right) +h\cdot
g\left( X,Y,v\right) .
\end{eqnarray*}%
\begin{equation}
\tag{5.5}
\end{equation}%
We now formulate a maximum principle for the optimal control system of (5.3).

\begin{theorem}
Suppose (A1)-(A2) hold. Let $\left( X\left( \cdot \right) ,Y\left( \cdot
\right) ,Z\left( \cdot \right) ,v^{*}\left( \cdot \right) \right) $ be an
optimal control and its corresponding trajectory of (5.3), $\left( p\left(
\cdot \right) ,q\left( \cdot \right) ,k\left( \cdot \right) ,h\left( \cdot
\right) \right) $ be the solution of (5.4). Then the maximum principle
holds, that is, for $t\in \left[ 0,T\right] $, $\forall v\in \mathcal{U},$ 
\begin{equation*}
\left\langle \bar H\left( t,X\left( t\right) ,Y\left( t\right) ,Z\left(
t\right) ,v^{*}\left( t\right) ,p\left( t\right) ,q\left( t\right) ,k\left(
t\right) ,h\left( t\right) \right) ,v-v^{*}\left( t\right) \right\rangle
\geq 0,\text{ a.e., a.s..}
\end{equation*}
\end{theorem}

\begin{proof}
By Theorem 6 in Section 3, we get the desired result.
\end{proof}

For relationship between (5.1) and (5.3), we have

\begin{lemma}
For any given admissible control $v\left( \cdot \right) ,$ we assume (A1)
and (A2) hold. Then (5.3) has a unique solution 
\begin{equation*}
\left( X^{t,x}\left( \cdot \right) ,Y^{t,x}\left( \cdot \right)
,Z^{t,x}\left( \cdot \right) \right) \in M^{2}\left( 0,T;\mathbb{R}^{n}%
\mathbb{\times R\times R}^{d}\right) .
\end{equation*}
\end{lemma}

\begin{lemma}
For any given admissible control $v\left( \cdot \right) ,$ we assume (A1)
and (A2) hold. Let 
\begin{equation*}
\left\{ u\left( t,x\right) ;0\leq t\leq T,x\in \mathbb{R}^{n}\right\} 
\end{equation*}%
be a random field such that $u\left( t,x\right) $ is $\mathcal{F}_{t,T}^{B}$%
-measurable for each $\left( t,x\right) ,$ $u\in C^{0,2}\left( \left[ 0,T%
\right] \times \mathbb{R}^{n};\mathbb{R}\right) $ a.s., and $u$ satisfies
SPDEs (5.1). Then $u\left( t,x\right) =Y^{t,x}\left( t\right) .$
\end{lemma}

\begin{lemma}
For any given admissible control $v\left( \cdot \right) ,$ we assume (A1)
and (A2) hold. Then 
\begin{equation*}
\left\{ u\left( t,x\right) =Y^{t,x}\left( t\right) ;0\leq t\leq T,x\in 
\mathbb{R}^{n}\right\} 
\end{equation*}%
is a unique classical solution of SPDEs (5.1).
\end{lemma}

The proofs are classical, we omit it. Now set the Hamilton function as
follows: 
\begin{eqnarray*}
\bar{H}\left( t,x,u,\nabla u\sigma ,v,p,q,k,h\right)  &=&l\left(
t,x,u,\nabla u\sigma ,v\right) -k\cdot G\left( t,x,u,\nabla u\sigma
,v\right)  \\
&&\ +q\cdot f\left( x,v\right) -p\cdot F\left( t,x,u,\nabla u\sigma
,v\right) +h\cdot g\left( x,v\right) .
\end{eqnarray*}%
We can state the maximum principle for the optimal control problem of SPDEs
(5.1).

\begin{theorem}[\textbf{Necessary maximum principle}]
Suppose $u\left( t,x\right) $ is the optimal solution of SPDEs (5.1)
corresponding to the optimal control $v^{\ast }\left( \cdot \right) $ of
(5.1). Then we have, for any $v\in \mathcal{U}$ and $t\in \left[ 0,T\right] ,
$ $x\in \mathbb{R}^{n}\mathbf{,}$ 
\begin{equation*}
\left\langle \bar{H}_{v}\left( t,x,u\left( t,x\right) ,\left( \nabla u\sigma
\right) \left( t,x\right) ,v^{\ast }\left( t\right) ,p\left( t\right)
,q\left( t\right) ,k\left( t\right) ,h\left( t\right) \right) ,v-v^{\ast
}\left( t\right) \right\rangle \geq 0,\text{ a.e., a.s.}
\end{equation*}
\end{theorem}

\begin{proof}
By virtue of lemma 9, 10, and 11, the optimal control problem of SPDE (5.1)
can be transformed into the one of FBDSDE (5.3). Hence, from Theorem 8, the
desired result is easily obtained.
\end{proof}

Next we apply our sufficient maximum principle to get the following result.

\begin{theorem}[\textbf{Sufficient maximum principle}]
For $\forall t\in \left[ 0,T\right] ,$ let $\hat{v}=\hat{v}\left( t\right)
\in \mathcal{U}_{ad}$ with corresponding solution $\hat{u}\left( t,x\right) $
of (5.1) and let $\left( \hat{X}\left( t\right) ,\hat{Y}\left( t\right) ,%
\hat{Z}\left( t\right) ,\hat{v}\left( t\right) \right) $ be quadruple and $%
\left( \hat{p}\left( t\right) ,\hat{q}\left( t\right) ,\hat{k}\left(
t\right) ,\hat{h}\left( t\right) \right) $ be a solution of the associated
adjoint FBDSDEs (5.3) and (5.4), respectively. Assume that 
\begin{equation*}
\bar{H}\left( t,X,Y,Z,v,\hat{p}\left( t\right) ,\hat{q}\left( t\right) ,\hat{%
k}\left( t\right) ,\hat{h}\left( t\right) \right)
\end{equation*}%
is convex in $\left( X,Y,Z,v\right) ,$ and $\gamma \left( Y\right) $ is
convex in $Y$, moreover the following condition holds 
\begin{eqnarray*}
&&\ \mathbb{E}\left[ \bar{H}\left( t,\hat{X}\left( t\right) ,\hat{Y}\left(
t\right) ,\hat{Z}\left( t\right) ,\hat{v}\left( t\right) ,\hat{p}\left(
t\right) ,\hat{q}\left( t\right) ,\hat{k}\left( t\right) ,\hat{h}\left(
t\right) \right) \right] \\
&=&\inf\limits_{v\in \mathcal{U}}\mathbb{E}\left[ \bar{H}\left( t,\hat{X}%
\left( t\right) ,\hat{Y}\left( t\right) ,\hat{Z}\left( t\right) ,v,\hat{p}%
\left( t\right) ,\hat{q}\left( t\right) ,\hat{k}\left( t\right) ,\hat{h}%
\left( t\right) \right) \right] .
\end{eqnarray*}%
Then $\hat{v}\left( t\right) $ is an optimal control for the problem (5.2).
\end{theorem}

\begin{proof}
Noting the above assumptions, by Theorem 7, it is easy to get desired result.
\end{proof}

\begin{remark}
In [13], Bernt Øksendal proved a sufficient maximum principle for the
optimal control of system described by a quasilinear stochastic heat
equation, that is 
\begin{eqnarray*}
dY\left( t,x\right)  &=&\left\{ 
\begin{array}{l}
\left[ LY\left( t,x\right) +b\left( t,x,Y\left( t,x\right) ,v\left( t\right)
\right) \right] dt \\ 
+\sigma \left( t,x,Y\left( t,x\right) ,v\left( t\right) \right) d%
\overrightarrow{W_{t}};%
\end{array}%
\right.  \\
\qquad \qquad \left( t,x\right)  &\in &\left[ 0,T\right] \times G.
\end{eqnarray*}%
\begin{equation}
\tag{5.6}
\end{equation}%
\begin{equation}
Y\left( 0,x\right) =\xi \left( x\right) ;\qquad x\in \overline{G}  \tag{5.7}
\end{equation}%
\begin{equation}
Y\left( t,x\right) =\eta \left( t,x\right) ;\qquad \left( t,x\right) \in
\left( 0,T\right) \times \partial G.  \tag{5.8}
\end{equation}%
Here $G$ is an open set in $\mathbb{R}^{n}$ with $C^{1}$ boundary $\partial G
$ and 
\begin{equation*}
L\phi \left( x\right) =\sum\limits_{i,j=1}^{n}a_{ij}\left( x\right) \frac{%
\partial ^{2}}{\partial x_{i}\partial x_{j}}\phi
+\sum\limits_{i=1}^{n}b_{i}\left( x\right) \frac{\partial }{\partial x_{i}}%
\phi ,\qquad \phi \in C^{2}\left( \mathbb{R}^{n}\right) 
\end{equation*}%
where $a\left( x\right) =\left[ a_{ij}\left( x\right) \right] _{1\leq
i,j\leq n}$ is a given symmetric definite symmetric $n\times n$ matrix with
entries $a_{ij}\left( x\right) \in C^{2}\left( G\right) \cap C\left( 
\overline{G}\right) $ for all $i,j=1,2,\cdots ,n$ and $b_{i}\left( x\right)
\in C^{2}\left( G\right) \cap C\left( \overline{G}\right) $ for all $%
i,j=1,2,\cdots ,n.$ It is worth to pointing out that our method to get the
sufficient maximum principle is completely different from his, and the most
important thing is that in our SPDEs, the coefficients of the elliptic
operator contain control variables (for more information see Theorem
2.1-Theorem 2.3 in [13]).
\end{remark}

\section{Applications}

We now illustrate the results of Section 3 by looking at some examples.
Theoretically, the maximum principles presented in Section 3 and Section 4
characterizes the optimal control through some necessary and sufficient
conditions. However, it is not immediately feasible to implement such
principles directly, partially due to the difficulty of computing fully
coupled forward-backward doubly stochastic system. In this section, we give
two special examples and show how to explicitly solve them using our maximum
principle.

\subsection{Example 1}

We provide a concrete example of forward-backward doubly stochastic LQ
problems and give the explicit optimal control and validate our major
theoretical results in Theorem 6. (Necessary maximum principle). First let
the control domain be $\mathcal{U}=\left[ -1,1\right] .$ Consider the
following linear forward-backward doubly stochastic control system. We
assume that $l=d=1.$%
\begin{equation}
\left\{ 
\begin{array}{l}
\text{\textrm{d}}y\left( t\right) =\left( z\left( t\right) -Z\left( t\right)
+v\left( t\right) \right) \text{\textrm{d}}\overrightarrow{W_{t}}-z\left(
t\right) \text{d}\overleftarrow{B_{t}}, \\ 
\text{\textrm{d}}Y\left( t\right) =-\left( z\left( t\right) +Z\left(
t\right) +v\left( t\right) \right) \text{\textrm{d}}\overleftarrow{B_{t}}%
+Z\left( t\right) \text{d}\overrightarrow{W_{t}}, \\ 
y\left( 0\right) =0,\quad Y\left( T\right) =0,\text{ }\quad t\in \left[ 0,T%
\right] ,%
\end{array}%
\right.  \tag{6.1}
\end{equation}%
where $T>0$ is a given constant and the cost function is 
\begin{eqnarray*}
J\left( v\left( \cdot \right) \right) &=&\frac{1}{2}\mathbb{E}%
\int_{0}^{T}\left( y^{2}\left( t\right) +Y^{2}\left( t\right) +z^{2}\left(
t\right) +Z^{2}\left( t\right) +v^{2}\left( t\right) \right) \text{d}t \\
&&+\frac{1}{2}\mathbb{E}Y^{2}\left( 0\right) +\frac{1}{2}\mathbb{E}%
y^{2}\left( T\right) .
\end{eqnarray*}%
\begin{equation}
\tag{6.2}
\end{equation}%
Note that (6.1) are linear control system. According to the existence and
uniqueness of (6.1), it is straightforward to know the optimal control is $%
u\left( \cdot \right) \equiv 0,$ with the optimal state trajectory $\left(
y\left( t\right) ,Y\left( t\right) ,z\left( t\right) ,Z\left( t\right)
\right) \equiv 0,$ $t\in \left[ 0,T\right] .$ Notice that the adjoint
equation associated with the optimal quadruple $\left( y\left( t\right)
,Y\left( t\right) ,z\left( t\right) ,Z\left( t\right) \right) \equiv 0$ are 
\begin{equation}
\left\{ 
\begin{array}{l}
\text{d}p\left( t\right) =-Y\left( t\right) \text{d}t+\left( -k\left(
t\right) -h\left( t\right) -Z\left( t\right) \right) \text{d}\overrightarrow{%
W_{t}}-k\left( t\right) \text{d}\overleftarrow{B_{t}}, \\ 
\text{d}q\left( t\right) =-y\left( t\right) \text{d}t+\left( -k\left(
t\right) -h\left( t\right) -z\left( t\right) \right) \text{d}\overleftarrow{%
B_{t}}+h\left( t\right) \text{d}\overrightarrow{W_{t}}, \\ 
p\left( 0\right) =0,\quad q\left( T\right) =0,\quad t\in \left[ 0,T\right] .%
\end{array}%
\right.  \tag{6.3}
\end{equation}%
Obviously, $\left( p\left( t\right) ,q\left( t\right) ,k\left( t\right)
,h\left( t\right) \right) \equiv 0$ is the unique solution of (6.3).
Instantly, we give the Hamiltonian function is 
\begin{eqnarray*}
&&H\left( t,y\left( t\right) ,Y\left( t\right) ,z\left( t\right) ,Z\left(
t\right) ,v,p\left( t\right) ,q\left( t\right) ,k\left( t\right) ,h\left(
t\right) \right) \\
&=&\frac{1}{2}\left( y^{2}\left( t\right) +Y^{2}\left( t\right) +z^{2}\left(
t\right) +Z^{2}\left( t\right) +v^{2}\right) \\
&&-k\left( t\right) \left( z\left( t\right) +Z\left( t\right) +v\right) \\
&&+h\left( t\right) \left( z\left( t\right) -Z\left( t\right) +v\right) \\
&=&\frac{1}{2}v^{2}.
\end{eqnarray*}%
It is clear that, for any $v\in \mathcal{U}$, we always have 
\begin{equation*}
\mathbb{E}\left\langle H_{v}\left( t,y\left( t\right) ,Y\left( t\right)
,z\left( t\right) ,Z\left( t\right) ,u\left( t\right) ,p\left( t\right)
,q\left( t\right) ,k\left( t\right) ,h\left( t\right) \right) ,v-u\left(
t\right) \right\rangle =0.
\end{equation*}

\subsection{Example 2}

In this subsection we will provide a special optimal control of SPDEs by
Theorem 13. (Sufficient maximum principle). We now introduce some notations.
For any random variable $F$ of the form 
\begin{equation*}
F=f\left( W\left( h_{1}\right) ,\ldots W\left( h_{n}\right) ;B\left(
k_{1}\right) ,\ldots B\left( k_{p}\right) \right)
\end{equation*}
with 
\begin{equation*}
f\in C_{b}^{\infty }\left( R^{n+p}\right) ,h_{1},\ldots h_{n}\in L^{2}\left( 
\left[ 0,T\right] ,R^{d}\right) ,k_{1},\ldots k_{p}\in L^{2}\left( \left[ 0,T%
\right] ,R^{l}\right) ,
\end{equation*}%
where 
\begin{equation*}
W\left( h_{i}\right) =\int_{0}^{T}h_{i}\left( t\right) \text{d}W_{t},\quad
B\left( h_{i}\right) =\int_{0}^{T}k_{i}\left( t\right) \text{d}B_{t},
\end{equation*}%
we let 
\begin{equation*}
D_{t}F=\sum_{i=1}^{n}f_{i}^{^{\prime }}\left( W\left( h_{1}\right) ,\ldots
W\left( h_{n}\right) ;B\left( k_{1}\right) ,\ldots B\left( k_{p}\right)
\right) h_{i}\left( t\right) ,\quad 0\leq t\leq T.
\end{equation*}%
For such an $F$, we define its 1,2-norm as: 
\begin{equation*}
\left\Vert F\right\Vert _{1,2}=\left( \mathbb{E}\left[ F^{2}+\int_{0}^{T}%
\left\vert D_{t}F\right\vert ^{2}\text{d}t\right] \right) ^{\frac{1}{2}}.
\end{equation*}%
$S$ denotes the set of random variable of the above form. We define the
Sobolev space: 
\begin{equation*}
\mathbb{D}^{1,2}=\overline{S}^{\left\Vert \cdot \right\Vert _{1,2}}.
\end{equation*}%
The \textquotedblright derivation operator\textquotedblright\ $D_{\cdot 
\text{ }}$extends as an operator from $\mathbb{D}^{1,2}$ into $L^{2}\left(
\Omega ;L^{2}\left( \left[ 0,T\right] ,R^{n}\right) \right) .$

Now we modify the stochastic reaction-diffusion equation considered in [13]
which can be described the density of a population at time $t\in \left[ 0,T%
\right] $ and at the point $x\in R$ as follows. 
\begin{equation}
\left\{ 
\begin{array}{c}
u\left( t,x\right) =x+\int_{t}^{T}\left[ v^{2}\left( s\right) \bigtriangleup
u\left( s,x\right) +u\left( s,x\right) +\nabla u\left( s,x\right) v\left(
s\right) \right] \text{d}s \\ 
+\int_{t}^{T}u\left( s,x\right) \text{d}\overleftarrow{B_{s}},\quad 0\leq
t\leq T,%
\end{array}%
\right.   \tag{6.4}
\end{equation}%
and $x\in R,$ $v\in \mathcal{U}_{ad}.$ The two Brownian motions $W$ and $B$
are one-dimensional. Suppose we want to minimize the following performance
criterion 
\begin{equation*}
J\left( v\right) =\mathbb{E}\left[ \int_{0}^{T}\frac{v^{\gamma }\left(
s\right) }{\gamma }\text{d}s+u\left( 0,x\right) \right] ,
\end{equation*}%
where $\gamma \geq 1.$ In this case the Hamiltonian gets the form 
\begin{eqnarray*}
&&\ H\left( t,X,Y,Z,v,p,q,k,h\right)  \\
\  &=&\frac{v^{\gamma }}{\gamma }-k\left( Y+Z\right) -pY+hv.
\end{eqnarray*}%
Obviously, it is convex in $\left( Y,Z,v\right) .$ The corresponding FBDSDEs
are 
\begin{equation}
\left\{ 
\begin{array}{l}
X^{t,x}\left( s\right) =x+\int_{t}^{s}v\left( r\right) \text{d}%
\overrightarrow{W_{r}}, \\ 
Y^{t,x}\left( s\right) =X^{t,x}\left( T\right) +\int_{s}^{T}\left(
Y^{t,x}\left( r\right) +Z^{t,x}\left( r\right) \right) \text{d}r \\ 
\qquad \qquad +\int_{s}^{T}Y\left( r\right) \text{d}\overleftarrow{B_{r}}%
-\int_{s}^{T}Z^{t,x}\left( r\right) \text{d}\overrightarrow{W_{r}},\quad
0\leq t\leq s\leq T,%
\end{array}%
\right.   \tag{6.5}
\end{equation}%
It is easy to obtain the solutions of (6.5) are 
\begin{equation}
Y^{t,x}\left( s\right) =\mathbb{E}\left[ \left. X^{t,x}\left( T\right) \exp
\left\{ W_{T}-W_{t}+B_{T}-B_{s}\right\} \right\vert \mathcal{F}_{s}\right] .
\tag{6.6}
\end{equation}%
Besides, the adjoint processes are 
\begin{equation}
\left\{ 
\begin{array}{l}
\text{d}p\left( s\right) =\left( p\left( s\right) +k\left( s\right) \right) 
\text{d}s+p\left( t\right) \text{d}\overrightarrow{W_{t}}-k\left( s\right) 
\text{d}\overleftarrow{B_{s}}, \\ 
\text{d}q\left( s\right) =h\left( s\right) \text{d}\overrightarrow{W_{s}},
\\ 
p\left( t\right) =-1,\quad q\left( T\right) =-p\left( T\right) ,\text{\qquad 
}t\leq s\leq T.%
\end{array}%
\right.   \tag{6.7}
\end{equation}%
The solutions of (6.7) are 
\begin{eqnarray*}
p\left( s\right)  &=&\mathbb{E}\left[ \left. -\exp \left\{
W_{s}+W_{t}+B_{s}-B_{t}\right\} \right\vert \mathcal{F}_{s}\right] , \\
q\left( s\right)  &=&\mathbb{E}\left[ \left. -p\left( T\right) \right\vert 
\mathcal{F}_{s}^{W}\right] , \\
h\left( s\right)  &=&D_{s}q\left( s\right) ,\quad \text{a.e., }0\leq t\leq
s\leq T.
\end{eqnarray*}%
\begin{equation}
\tag{6.8}
\end{equation}%
The function 
\begin{eqnarray*}
v &\rightarrow &H\left( t,X,Y,Z,v,p,q,k,h\right)  \\
\  &=&\frac{v^{\gamma }}{\gamma }-kY-pY+hv.
\end{eqnarray*}%
is minimum when 
\begin{equation*}
v\left( t\right) =\left( h\left( t\right) \right) ^{\frac{1}{\gamma -1}%
},\quad 0\leq t\leq T.
\end{equation*}%
where $h\left( t\right) $ are given by (6.8).


\begin{thebibliography}{10}
\bibitem[1]{Ben} A. Bensoussan, Lectures on Stochastic Control, Lecture
Notes in Mathematics, Vol. 972, Nonlinear Filtering and Stochastic Control,
Proceeding, Cortona, 1981.

\bibitem[2]{Ben} A. Bensoussan, Stochastic maximum principle for distributed
parameter system. J. Franklin Inst\textit{.} 315\textbf{\ }(1983),. 387--406.

\bibitem[3]{Ben} A. Bensoussan, Stochastic Control of Partially Observable
Systems. Cambridge Uni- versity Press 1992.

\bibitem[4]{Bis} J. M. Bismut, An introductory approach to duality in
optimal stochastic control. SIAM Rev\textit{.,} 20 (1978),. 62--78.

\bibitem[5]{Has} U. G. Haussmann,\ General necessary conditions for optimal
control of stochastic system, Maht. Programm. Stud., 6 (1976), 34--48.

\bibitem[6]{Has} U. G. Haussmann, A stochastic maximum principle for optimal
control of diffusions. Pitman Research Notes in Mathematics 151 (1987).

\bibitem[7]{HSW} Y. Han, S. Peng, Z. Wu, Maximum Principle for Backward
Doubly Stochastic Control Systems with Applications, SIAM J. Control Optim.,
48(7), 4224-4241.

\bibitem[8]{JZ} S. Ji and X.Y. Zhou, A maximum principle for stochastic
optimal control with terminal state constraints, and its applications.
Communications in Information and Systems 6 (4) (2006) 321--338.

\bibitem[9]{K} H.J. Kushner, Necessary conditions for continuous parameter
stochastic optimization problems, SIAM J. Control, 10 (1972), 550--565.

\bibitem[10]{M} R. E. Mortensen, Stochastic optimal control with noisy
observations. Int. J. Control 4 (1966), 455--464.

\bibitem[11]{NP} D. Nualart and E. Pardoux, Stochastic calculus with
anticipating integrands, Probab. Theory Related Fields, 78 (1988), 535--581.

\bibitem[12]{N} M. Nisio, Optimal control for stochastic partial
differential equations and viscosity solutions of bellman equations Nagoya
Math. J. Vol. 123 (1991), 13-37.

\bibitem[13]{O} B. Øksendal, Optimal Control of Stochastic Partial
Differential Equations. Stochastic Anal. Appl. \textbf{23 }(2005), 165--179.

\bibitem[14]{PP1} E. Pardoux and S. Peng, Adapted solution of a backward
stochastic differential equation.\textit{\ }Systems Control Letters. 14
(1990). 55--61.

\bibitem[15]{PP2} E. Pardoux and S. Peng, Backward doubly stochastic
differential equations and systems of quasilinear parabolic SPDEs, Probab.
Theory Related Fields, 98 (1994), 209--227.

\bibitem[16]{P1} S. Peng, A general stochastic maximum principle for optimal
control problems, SIAM J. Control, 28 (1990), 966--979.

\bibitem[17]{P2} S. Peng, Backward stochastic differential equations and
application to optimal control. Applied Mathematics and Optimization 27 (4)
(1993), 125-144.

\bibitem[18]{PS} S. Peng and Y. Shi, A Type of Time-Symmetric
Forward-Backward Stochastic Differential Equations. C. R. Acad. Sci. Paris,
Ser. I 336 (9) (2003), 773-778.

\bibitem[19]{PW} S. Peng and Z. Wu, Fully Coupled Forward-Backward
Stochastic Differential Equations and Applications to Optimal Control\textit{%
. }SIAM J. Control Optim. 37 (1999), 825-843.

\bibitem[20]{Pon} L.S. Pontryagin, V.G. Boltyanskti, R.V. Gamkrelidze, E.F.
Mischenko, The Mathematical Theory of Optimal Control Processes.
Interscience, John Wiley, New York (1962).

\bibitem[21]{S} J. Shi and Z. Wu, The maximum principle for fully coupled
forward-backward stochastic control system. Acta Automatica Sinica 32 (2)
(2006), 161-169.

\bibitem[22]{W} Z. Wu, Maximum principle for optimal control problem of
fully coupled forward-backward stochastic systems. Systems Sci. Math. Sci.
11 (3) (1998), 249-259.

\bibitem[23]{X} W. Xu, Stochastic maximum principle for optimal control
problem of forward and backward system. J. Australian Mathematical Society
B37 (1995), 172-185.

\bibitem[24]{YZ} J. Yong and X.Y. Zhou, Stochastic Controls: Hamiltonian
Systems and HJB Equations. Springer, New York 1999.

\bibitem[25]{ZZ} Q. Zhang and H. Zhao, Stationary solutions of SPDEs and
infinite horizon BDSDEs, J. Funct. Anal., 252 (2007), 171-219.

\bibitem[26]{ZS} L. Zhang and Y. Shi, Maximum Principle for Forward-Backward
Doubly Stochastic Control Systems and Applications. ESAIM: COCV.
DOI:10.1051/cocv/2010042.
\end{thebibliography}
\end{document}